\documentclass[12pt]{article}

\usepackage[utf8]{inputenc}
\usepackage{epsfig}
\usepackage{amsmath}
\usepackage{amsfonts}
\usepackage{amsthm}

\usepackage{algpseudocode}
\usepackage{algorithm}

\usepackage{color}

        \theoremstyle{plain}

        \newtheorem{prop}{Property}[section]

        \theoremstyle{remark}
        \newtheorem{remark}{\bf Remark}[section]
        \theoremstyle{remark}
        
        \theoremstyle{remark}
        
        \newtheorem{algo}{Algorithm}[section]

\newcommand{\demi}{\frac{1}{2}}

\newcommand{\R}{{\mathbb{R}}}

\newcommand{\dt}{\partial_t}
\newcommand{\dx}{\partial_x}

\newcommand{\cint}[1]{\left\langle #1\right\rangle}
\newcommand{\x}{{\bf x}}

\newcommand{\vit}{{\bf v}}
\newcommand{\um}{{\bf u}}

\newcommand{\Maxw}{M}

\newcommand{\Dv}{\Delta v}

\newcommand{\Dt}{\Delta t}

\newcommand{\fnijkq}{f^n_{ijkq}}
\newcommand{\fnpunijkq}{f^{n+1}_{ijkq}}
\newcommand{\gnijkq}{g^n_{ijkq}}
\newcommand{\gnpunijkq}{g^{n+1}_{ijkq}}

\newcommand{\Unijk}{U^n_{ijk}}
\newcommand{\Unpunijk}{U^{n+1}_{ijk}}

\newcommand{\modif}[1]{#1}

\begin{document}


\begin{center}

{\bf Locally refined
  discrete velocity grids for stationary rarefied flow simulations }

\vspace{1cm}
C. Baranger$^1$, J. Claudel$^1$, N. H\'erouard$^{1,2}$, L. Mieussens$^2$

\bigskip
$^1$CEA-CESTA\\
15 avenue des sabli\`eres - CS 60001\\
33116 Le Barp Cedex, France\\
{ \tt(celine.baranger@cea.fr, jean.claudel@cea.fr, nicolas.herouard@cea.fr)}\\

\bigskip
$^2$Univ. Bordeaux, IMB, UMR 5251, F-33400 Talence, France.\\
CNRS, IMB, UMR 5251, F-33400 Talence, France. \\
INRIA, F-33400 Talence, France. \\

{ \tt(Luc.Mieussens@math.u-bordeaux1.fr)}

\end{center}

Abstract: Most of deterministic solvers for rarefied gas dynamics use
discrete velocity (or discrete ordinate) approximations of the
distribution function on a Cartesian grid. This grid must be
sufficiently large and fine to describe the distribution functions at
every space position in the computational domain. For 3-dimensional
hypersonic flows, like in re-entry problems, this induces much too
dense velocity grids that cannot be practically used, for memory
storage requirements. In this article, we present an approach to
generate automatically a locally refined velocity grid adapted to a
given simulation. This grid contains much less points than a standard
Cartesian grid and allows us to make realistic 3-dimensional
simulations at a reduced cost, with a comparable accuracy. 


\bigskip


\section{Introduction}

The description of a flow surrounding a flying object at hypersonic
speed in the rarefied atmosphere (more than 60 km altitude) is still a
challenge in the atmospheric Re-Entry
community~\cite{anderson}. When
this flow is in a rarefied state, that is to say when the Knudsen
number (which is the ratio $Kn=\frac{\lambda}{L}$ between the mean
free path $\lambda$ of particle and a characteristic macroscopic
length $L$) is larger than $0.01$, the flow cannot be accurately
described by the compressible Navier-Stokes equations of Gas
Dynamics. In this case, the kinetic theory has to be used. The
evolution of the molecules of the gas is then described by a mass
density distribution in phase space, which is a solution of the
Boltzmann equation. In the transitional regime, this equation can be
replaced by the simpler Bhatnagar-Gross-Krook (BGK) model.

In order to be able to compute parietal heat flux and aerodynamic
coefficients in the range of 60-120 km, a
kinetic description of the stationary flow is necessary. 

The most popular numerical method to simulate rarefied flows is the
Direct Simulation Monte Carlo method (DSMC)~\cite{bird}. However, it
is well known that this method is very expensive in transitional
regimes, in particular for flows in the range of altitude we are
interested in. The efficiency of DSMC can be improved by using
coupling strategies (see~\cite{BB2009,DD2012}) or implicit schemes
(see~\cite{PR_2001,DP2008}), but these methods are still not very well
suited for stationary computations. In contrast, deterministic methods
(based on a numerical discretization of the stationary kinetic model)
can be more efficient in transitional regimes. Up to our knowledge,
there are few deterministic simulation codes specifically designed for
steady flows. One of the most advanced ones is the 3D code of
Titarev~\cite{T_cicp} developed for unstructured meshes. Another 3D
code has been developed by G. Brook~\cite{PhD_Brook}. Other codes
exist, but they are rather designed for unsteady problems, see for
instance~\cite{KAAFZ,AC2011} or the recent UGKS scheme developed by
K. Xu and his collaborators~\cite{XH2010,HXY2012}, which is an
Asymptotic Preserving scheme for unsteady flows.

In our team, we developed several years ago a code to make 2D plane
and axisymmetric simulations of rarefied flows based on the BGK model
(see a description in~\cite{luc_m3as,luc_jcp,ADM2009}). This code has
recently been extended to 3D computations, for polyatomic gases.
 Due to the physical model (polyatomic gases), the space discretization (block structured
mesh), and the parallelization (space domain decomposition with MPI
and inner parallelization with openMP), this code is rather different
from the other existing 3D codes recently presented in the literature
for the same kind of problems (the 3D code of Titarev~\cite{T_cicp}
for example), even if space domain decompositions have
  already been used for unsteady simulation (see \cite{KAAFZ}).

  All the codes designed for steady flows have a common feature: they
  are based on a ``discrete ordinate'' like approach, and use a global
  velocity grid. This grid is generally a Cartesian grid with a
  constant step size. The number of points of this grid is roughly
  proportional to the Mach number of the flow in each direction, and
  hence can be prohibitively large for hypersonic flows, even with
  parallel computers. To compute realistic cases (3D configurations
  with Mach number larger than 20), the velocity space discretization
  has to be modified in order to decrease CPU time and memory storage
  requirements.  It has already been noticed that a refinement of the
  grid around small velocities can improve the accuracy and reduce the
  cost of the computation (especially for large Knudsen numbers in
  flows close to solid boundaries, see~\cite{T_cicp}).  However, up to
  our knowledge, there is no general strategy in the literature that
  helps us to reduce the number of discrete velocities of a velocity
  grid for any rarefied steady flow, even if some works on adaptive
  velocity grids have already been presented: the first attempt seems
  to be~\cite{Aristov77} for a 1D shock wave problem,
  and recently, more general adaptive grid techniques designed for
  unsteady flows have been presented
  in~\cite{KAF2011,CXLC_2012,KA2012,BM2013}.

  The main contribution of this article is to propose an algorithm for
  an automatic construction of a locally refined velocity grid that
  allows a dramatic reduction of the number of discrete velocities,
  with the same accuracy as a standard Cartesian grid. This algorithm
  uses a compressible Navier-Stokes pre-simulation to obtain a rough
  estimation of the macroscopic fields of the flow. These fields are
  used to refine the grid wherever it is necessary by using some
  indicators of the width of the distribution functions in all the
  computational domain. This strategy allows us to simulate hypersonic
  flows that can hardly be simulated by standard methods, since we are
  indeed able to apply our method to our kinetic code to simulate a
  re-entry flow at Mach 25 and for temperature and pressure conditions
  of an altitude of 90 km. In this example, the CPU time and memory
  storage can be decreased up to a factor 24, as compared to a method
  with a standard Cartesian velocity grid. Note that preliminary results have
  already been presented in~\cite{BCHM_vancouver}
  and~\cite{BCHM_saragosse}.

The outline of this article is the following. In
section~\ref{sec:btz}, we briefly present the kinetic description of a
rarefied gas. In section~\ref{sec:algo}, we discuss the problems
induced by the use of a global velocity grid, and our algorithm is
presented. Our approach is illustrated in section~\ref{sec:res} with
several numerical tests. To simplify the reading of the paper, the
presentation of our simulation code is
made in the appendix. 


\section{Boltzmann equation and Cartesian velocity grid}
\label{sec:btz}

   \subsection{Kinetic description of rarefied gases}

In kinetic theory, a monoatomic gas is described by the distribution
function $f(t,\x,\vit)$ defined such that $f(t,\x,\vit)d\x d\vit$ is the mass
of molecules that at time $t$ are located in an elementary space volume
$d\x$ centered in $\x=(x,y,z)$ and have a velocity in an elementary
volume $d\vit$ centered in $\vit=(v_x,v_y,v_z)$.

Consequently, the macroscopic quantities as mass density $\rho$,
momentum $\rho \um$ and total energy $E$ are defined as the five first
moments of $f$ with respect to the velocity variable, namely:
\begin{equation}  \label{eq-mts}
(\rho(t,\x),\rho\um(t,\x),E(t,\x))=\int_{\R^3}(1,\vit,\demi|\vit|^2)f(t,\x,\vit)
\, d\vit.
\end{equation}

The temperature $T$ of the gas is defined by the relation
$E=\demi\rho|\um|^2+\frac{3}{2}\rho R T$, where $R$ is the gas constant
defined as the ratio between the Boltzmann constant and the molecular
mass of the gas.

When the gas is in a thermodynamical equilibrium state, it is well known that
the distribution function $f$ is a Gaussian function
$\Maxw[\rho,\um,T]$ of $\vit$, called
Maxwellian distribution, that depends only on the
macroscopic quantities :
\begin{equation}  \label{eq-maxw}
\Maxw[\rho,\um,T]=\frac{\rho}{(2\pi RT)^\frac{3}{2}}\exp(-\frac{|\vit-\um|^2}{2RT}).
\end{equation}
It can easily be checked that $\Maxw$ satisfies
relations~(\ref{eq-mts}).

If the gas is not in a thermodynamical equilibrium state, its
evolution is described by the following kinetic equation
\begin{equation}  \label{eq-kin}
\dt f + \vit \cdot \nabla_\x f = Q(f),
\end{equation}
which means that the total variation of $f$ (described by the
  left-hand side) is due to the collisions between molecules
  (described by the right-hand side). The most realistic collision
model is the Boltzmann operator but it is still very computationally
expensive to use. In this paper, we use the simpler BGK
model~\cite{bgk,welander}
\begin{equation}  \label{eq-bgk}
Q(f) = \frac{1}{\tau}(\Maxw[\rho,\um,T]-f)
\end{equation}
which models the effect of the collisions as a relaxation of $f$
towards the local equilibrium corresponding to the macroscopic
quantities defined by~(\ref{eq-mts}).
The relaxation time is defined
as $\tau=\frac{\mu}{\rho RT}$, where $\mu$ is the viscosity of the
gas. This definition allows to match the correct viscosity in the
Navier-Stokes equations obtained by the Chapman-Enskog~expansion. This
viscosity $\mu$ is usually supposed to fit the law
$\mu=\mu_{ref}(\frac{T}{T_{ref}})^\omega$, where $\mu_{ref}$ and
$T_{ref}$ are reference viscosity and temperature determined
experimentally for each gas, as well as the exponent $\omega$ (see a
table in~\cite{bird}). 

The interactions of the gas with solid boundaries are described with
\modif{the diffuse reflection model.}
Let us suppose that the boundary has a
velocity $\um_w=0$ and temperature $T_w$. In the diffuse
  reflection model, a molecule that
collides with this boundary is supposed to be re-emitted  with a
temperature equal to $T_w$, and with a random velocity normally
distributed around 0. This reads
\begin{equation}\label{eq-BCdiff} 
  f(t,\x,\vit) =  \Maxw[\sigma_w,0,T_w](\vit)
\end{equation}
if $\vit\cdot {\bf n}(\x) > 0$, where ${\bf n}(\x)$ is the normal to
the wall at point $\x$ directed into the gas.
The parameter $\sigma_w$
is defined so that there is no normal mass flux across the boundary (all
the molecules are re-emitted). Namely, that is,  
\begin{equation}  \label{eq-sigmadiff}
\sigma_w=-\left(\int_{\vit\cdot {\bf
    n}(\x) < 0}f(t,\x,\vit) \vit\cdot {\bf n}(\x) \, d\vit\right)/\left(\int_{\vit\cdot {\bf
    n}(\x) > 0}\Maxw[1,0,T_w](\vit) \vit\cdot {\bf n}(\x) \, d\vit\right).
\end{equation}
\modif{There are other reflection models, like the Maxwell model with partial
accomodation, but they are not used in this work.}








   \subsection{Velocity discretization on Cartesian grid}

\begin{figure}
  \centering\includegraphics[height=0.3\textheight]{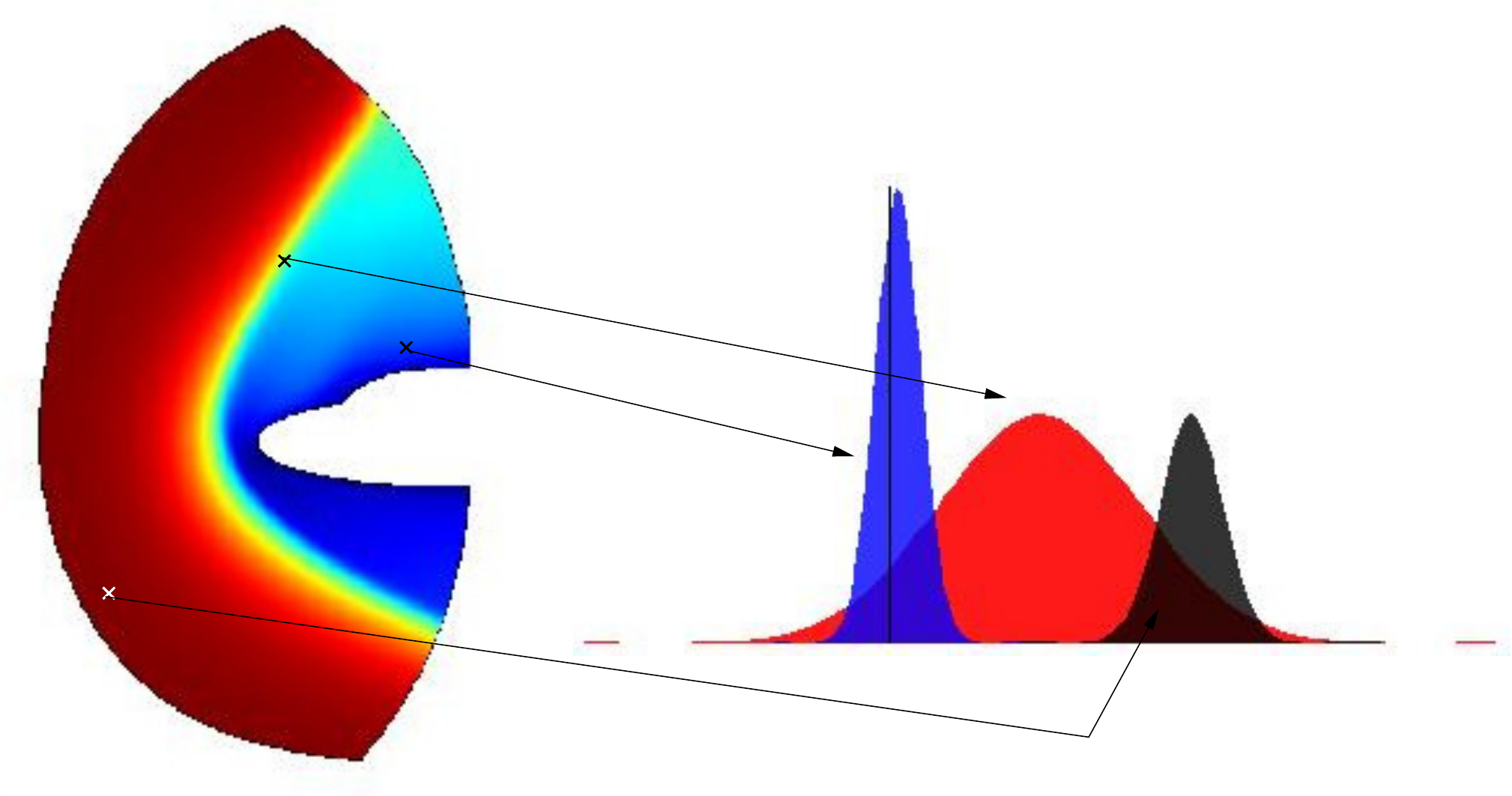}    \caption{Three distribution functions in different space points of a computational domain for a re-entry problem}
  \label{fig:3f}
\end{figure}

   In most of existing computational codes for the Boltzmann equation
   that use a deterministic method, the velocity variable is
   discretized with a {\it global} Cartesian grid, that is to say the
   same grid for every point of the space mesh. Consequently, it is
   necessary to compute a priori a velocity grid in which every
   distribution function that may appear in the computation can be
   resolved (see figure~\ref{fig:3f}). This means that the grid must
   be:
  \begin{itemize}
  \item large enough to contain the main part of the
    distribution functions {\it for every position in the
      computational domain}

  \item fine enough to capture the core of the distribution functions
    {\it for every position in the computational domain}
\end{itemize}

\begin{figure}
  \centering
  \includegraphics[height=0.3\textheight]{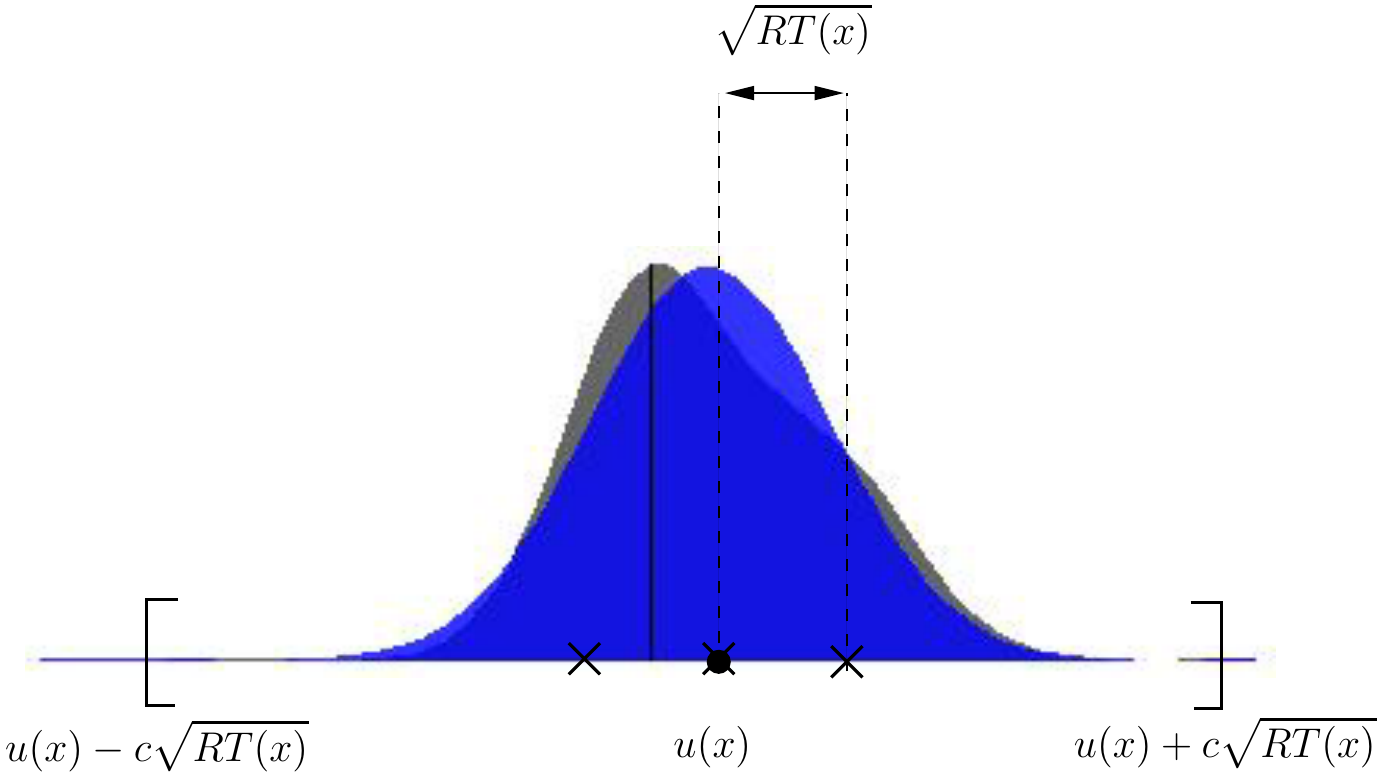} 
  \caption{Support of a distribution function $f(x,.)$: the distribution (in
    black), its corresponding Maxwellian (in blue), and the
    corresponding support based on the bulk velocity $u$ and the
    temperature $T$.}
  \label{fig:grille_param}
\end{figure}

\begin{figure}
  \centering
     \includegraphics[height=0.45\textheight]{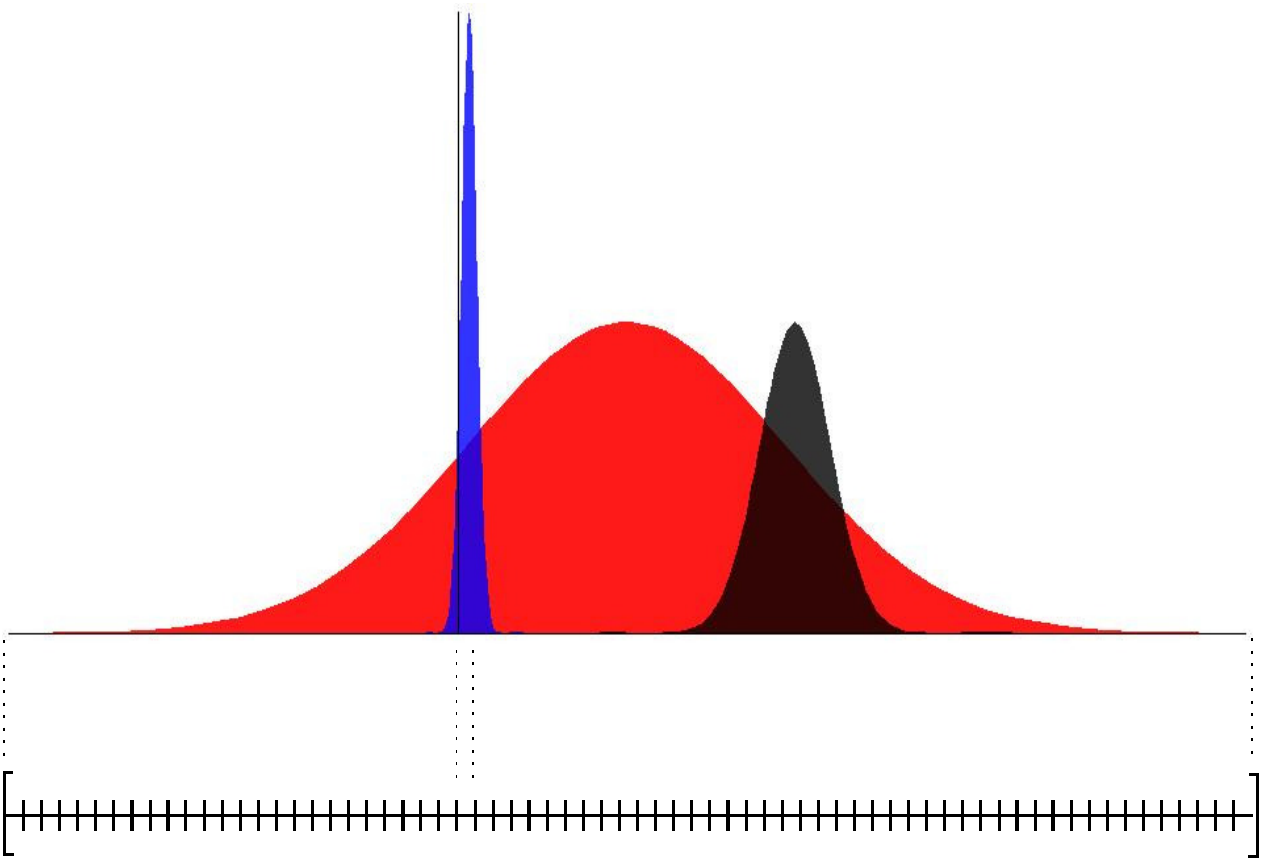}      
  \caption{Three different distribution functions and a global grid
    that satisfies~(\ref{eq-bounds})-(\ref{eq-step}): the bounds are
    given by the largest distribution, the step is given by the
    narrowest distribution.}
  \label{fig:grille_globale}
\end{figure}

The corresponding parameters (bounds and step of the grid) can be determined with the following idea: at each
point $\x$ of the computational domain, the macroscopic velocity $\um(\x)$
and temperature $T(\x)$ give information on the support of the
distribution function $f(\x,.)$. Indeed, if $f(\x,.)$ is ``not too far"
from its corresponding Maxwellian, as supposed in the BGK model, its support is centered around $\um(\x)$,
and is essentially localized between bounds that depend on $\um(\x)$
and $T(\x)$, corresponding to the standard deviation of the Maxwellian ; that is to say $f(\x,.)$ is very small outside (see
figure~\ref{fig:grille_param} where these bounds are
$\um(\x)-c\sqrt{RT(\x)}$ and $\um(\x)+c\sqrt{RT(\x)}$ in 1D). When several distribution functions
have to be discretized on the same grid, their supports are reasonably well approximated if the bounds are
\begin{equation}\label{eq-bounds} 
\begin{split}  
v^{\alpha}_{max}=\max_{x\in\Omega}\lbrace
u^{\alpha}(\x)+c\sqrt{RT(\x)} \rbrace \qquad  \qquad 
v^{\alpha}_{min}=\min_{\x\in\Omega}\lbrace
u^{\alpha}(\x)-c\sqrt{RT(\x)} \rbrace , 
\end{split}  
\end{equation}
and if the grid step is
\begin{equation}  \label{eq-step}
\Delta v = a\min_{\x\in\Omega}{\sqrt{RT(\x)}},
\end{equation} 
where $\Omega$ is the computational domain (in space), and
$\alpha=x,y,z$ (see~\ref{fig:grille_globale}). In~(\ref{eq-bounds})-(\ref{eq-step}), the parameters $c$ and $a$ can be
chosen as follows. For $c$, 
statistical argument suggest that values between 
$3$ and $4$ are needed, and so $c=3$ seems to be a good compromise between accuracy and computational
cost.
The parameter $a$ allows us to adjust the grid step: it must be
at most equal to 
1
(which ensures that there are at least 
six
points
inside the core of every distribution function), but a smaller value might be necessary for an
accurate computation.

To compute these bounds, it is necessary to a priori estimate the
values of $\um$ and $T$ in the computational domain. A simple way is
to use a compressible Navier-Stokes (CNS) pre-simulation (for which the CPU
time is negligible as compared to the Boltzmann simulation). This
simulation gives the fields $\um^{CNS}$ and $T^{CNS}$ in $\Omega$, and
the bounds and the grid step can be defined by
formulas~(\ref{eq-bounds})-(\ref{eq-step}) applied to these fields, that
is to say
\begin{equation}\label{eq-bounds_CNS} 
\begin{split}  
v^{\alpha}_{max}=\max_{\x\in\Omega}\lbrace u^{\alpha,CNS}(\x)+c\sqrt{RT^{CNS}(\x)}\rbrace \qquad  
v^{\alpha}_{min}=\min_{\x\in\Omega}\lbrace u^{\alpha,CNS}(\x)-c\sqrt{RT^{CNS}(\x)}\rbrace, 
\end{split}  
\end{equation}
and 
\begin{equation}  \label{eq-step_CNS}
\Delta v = a\min_{\x\in\Omega}{\sqrt{RT^{CNS}(\x)}}.
\end{equation} 

However, the values of the bounds are mainly determined by the temperatures in the shock zone, which can reach thousands of Kelvin.
One must be careful that, at high altitude and high speeds, numerical computations of Navier-Stokes equations can underestimate those temperatures,
which would lead to inappropriate bounds.

It is quite clear that the size of the velocity grid increases with
the Mach number. Indeed, a large Mach number implies large upstream
velocities and large temperature in the shock wave, which lead to very
large bounds (see~(\ref{eq-bounds})), while the temperature of the body
remains small, thus leading to small grid
step~(see~(\ref{eq-step})). For realistic re-entry problems, this can
lead to prohibitively large grids. For instance, for the flow
around a cylinder of radius 0.1 m at Mach 20 and altitude 90 km,
formulas~(\ref{eq-bounds_CNS}) and~(\ref{eq-step_CNS}) used with a Compressible
Navier-Stokes pre-simulation lead to a velocity grid of $52 \times 41
    \times 41$ points and around 350 GB memory requirements with a
    coarse 3D mesh in space. 

However, this is mainly due to the use of a Cartesian grid with a
constant step. In order to decrease the size of the grid, it is
attractive to refine it only wherever it is necessary, and to coarsen
it elsewhere. In the following section, we show that this can be done
automatically with our new algorithm.

This approach may not be well designed for cases where $f$ is very far from it corresponding Maxwellian (like in shock interactions problem). However, for reentry problems, our assumption is realistic for transitional regime we are interested in.

\section{An algorithm to define locally refined discrete velocity grids}
\label{sec:algo}

\begin{figure}
  \centering
     \includegraphics[height=0.45\textheight]{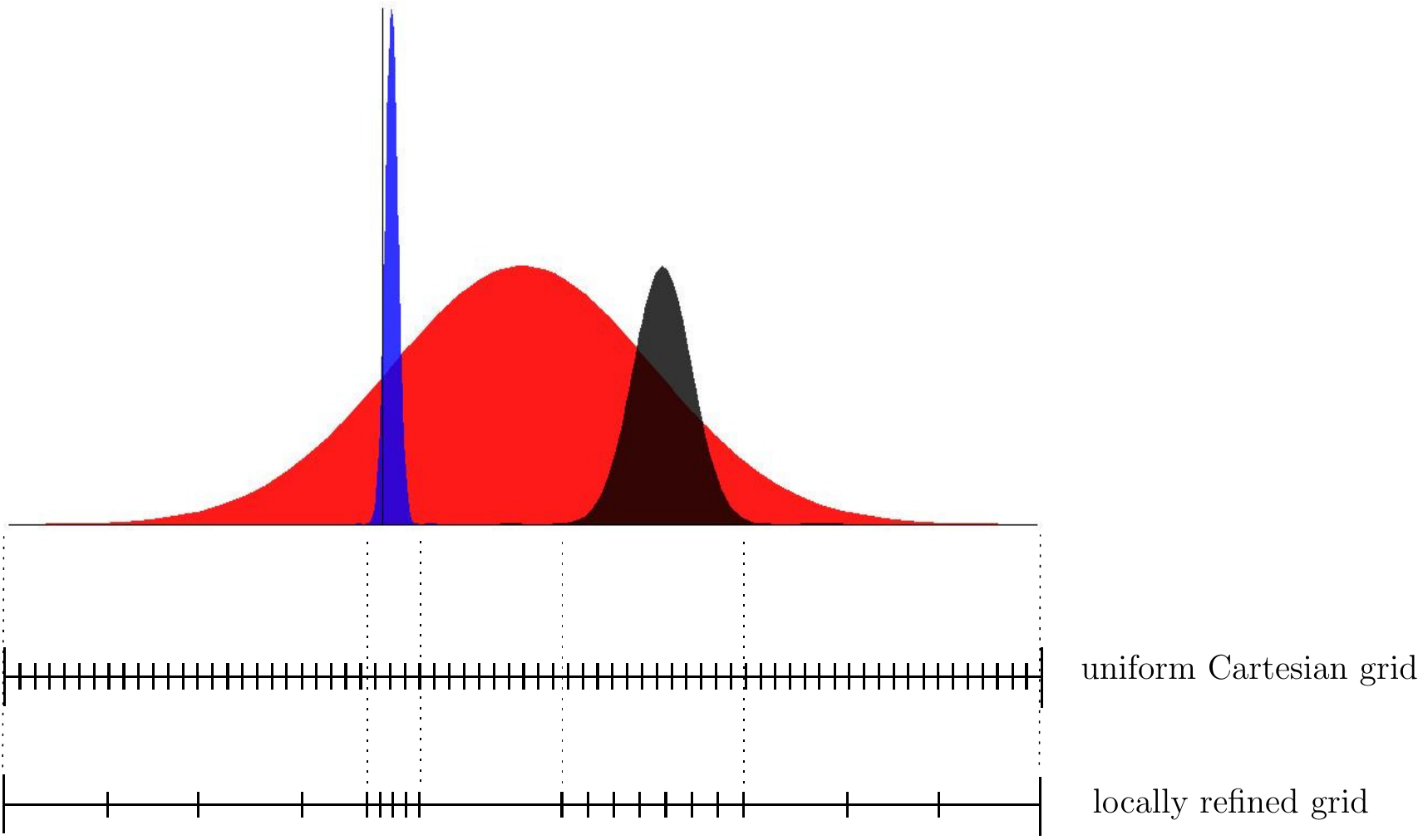}      
  \caption{Three different distribution functions (top), the
    corresponding uniform Cartesian velocity grid (middle), and the
    grid locally refined in the support of the
    distributions (coarsened elsewhere).}
  \label{fig:grille_raf}
\end{figure}

Since we have to represent many distribution functions on the same
grid, it is natural to refine the grid in the cores of theses
distributions, and to coarsen it in the tails (see
figure~\ref{fig:grille_raf}). To achieve this goal, we define the concept of
``support function'', and we use it to design an AMR (Adaptive Mesh
Refinement) velocity grid.

   \subsection{The support function}
\label{subsec:algo_support}

\begin{figure}
  \centering
\includegraphics[height=0.3\textheight]{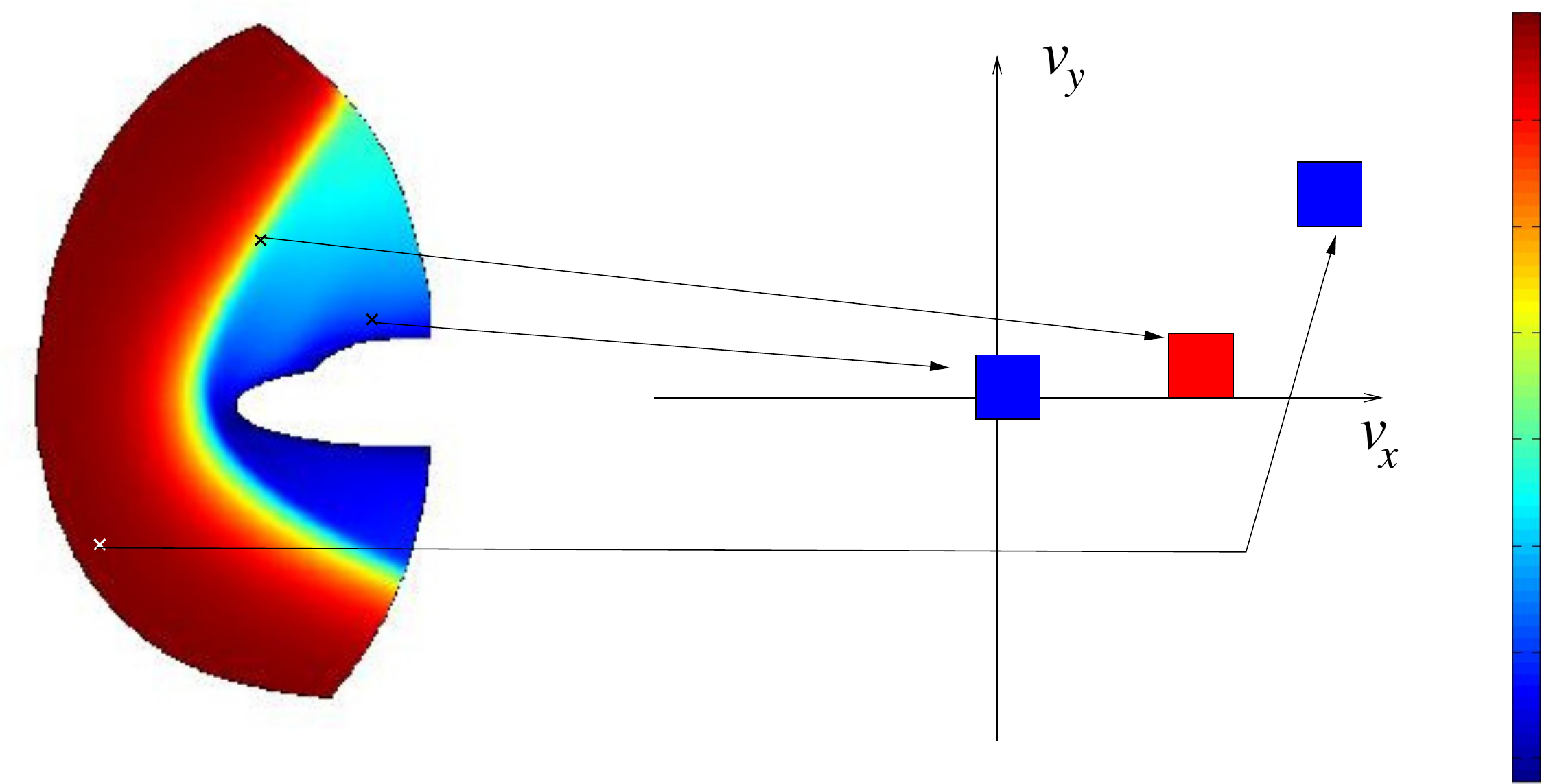} \\
 \caption{The macroscopic CNS flow (left), and some values of $\phi$
   at three different $v$ (right).}
  \label{fig:support}
\end{figure}

At each point $\x$ of the computational space domain $\Omega$, $\um(\x)$ and $T(\x)$ give
information on the ``main'' support of the distribution function $f(\x,.)$
in the velocity space: this support is centered at $\um(\x)$ and of
size $2c\sqrt{RT(\x)}$ (see~(\ref{eq-bounds})). We define the {\it
  support function} $\phi$ in the velocity space as follows: for
every $\vit$ in the velocity space, we set $\phi(\vit)=\sqrt{RT(\x)}$ if
there exists in $\Omega$ a $\x$ such that $||\vit-\um(\x)|| \leq c\sqrt{RT(\x)}$, that is to say that
 $\vit$ is inside the support of a distribution, considered as a sphere of center $\um(\x)$ and of radius $c\sqrt{RT(\x)}$. This function gives
an incomplete mapping of the velocity domain (see figure~\ref{fig:support}).
Now we want to extend this function to the whole velocity space, so
that it can be used as a refinement criterion to design an optimal
velocity grid.

Indeed, so far, this function is not defined for every $\vit$,
especially for large ones, as we cannot find for every $\vit$ an appropriate pair $(\um(\x),T(\x))$ to match the definition
of the support function. Moreover, this function might be multivalued, as
there might be two different space points $\x_1$ and $\x_2$ with same
velocity $\um(\x_1)=\um(\x_2)$, but with different temperatures
$T(\x_1)\ne T(\x_2)$. Since our goal is to obtain for every $\vit$ the
minimum size of the support of every distribution centered around $\vit$,
these two problems can be avoided as follows :
\begin{itemize}
  \item[(a)] we set $\phi(\vit)=\min(\sqrt{RT(\x_1)},\sqrt{RT(\x_2)})$ if
    $||\vit-\um(\x_1)|| \leq c\sqrt{RT(\x_1)}$ and $||\vit-\um(\x_2)||
    \leq c\sqrt{RT(\x_2)}$ with $(\um(\x_1),T(\x_1))\ne
    (\um(\x_2),T(\x_2))$. Indeed, the size of
    the grid around $\vit$ is constrained by the distribution centered on
    $\vit$ that has the smallest support, hence $\phi$ must have the
    corresponding smallest possible value.
  \item[(b)] 
 if there is no $\x$ such
    that $||\vit-\um(\x)|| \leq c\sqrt{RT(\x)}$, then $\vit$ is in the tail of
    every distribution function, and there is no reason to refine the
    grid there, so we set $\phi(\vit)$ to its largest possible value
    $\phi(\vit) = \max{\sqrt{RT(\x)}} $.

\end{itemize}
These ideas lead to the following algorithm for an automatic
construction of $\phi$. Note that below, the computational
domain in space is supposed to be discretized with a mesh composed of cells numbered
with three indices $(i,j,k)$. This is not a restriction of
our algorithm. 
\vspace{0.5cm}
\begin{algo}[Construction of the support function]
\label{algo:support}
$ \ $

  \begin{enumerate}
  \item CNS velocity and temperature are stored in arrays $\um(i,j,k)$ and $T(i,j,k)$ for $i,j,k=1:i_{max},j_{max},k_{max}$.
  \item construction of a (fine) Cartesian velocity grid (based on the
    previous global criterion~(\ref{eq-bounds_CNS})-(\ref{eq-step_CNS})): points
    $\vit(q)$ with $q=0:q_{max}$.
\item computation of the field $\sqrt{RT}$, stored in an descending
  order in the 1D array $\psi(I)$, for $I=1:i_{max}\times j_{max}\times
  k_{max}$, so that $\psi(I)=\sqrt{RT((i,j,k)(I))}$.
\item initialization of the array $\phi(0:q_{max})=\max(\psi)$ on the
  velocity grid (one value per discrete velocity)
\item  $do \ I=1:i_{max}\times  j_{max}\times k_{max}$ (loop on the cells of the space mesh)
    \begin{itemize}
       \item[]  $do \ q=0:q_{max}$ (loop on the nodes of the velocity grid)
           \begin{itemize}
              \item if $\vit(q)$ is in the sphere of center $\um((i,j,k)(I))$
                and radius $c\psi(I)$, then it is in the support of
                a distribution function, and we set $\phi(q):=\psi(I)$
           \end{itemize}

     \end{itemize}

\end{enumerate}
\end{algo}

This algorithm ensures that the array $\phi$ satisfies the following
property.
\begin{prop}
  For every $q$ between $0$ and $q_{max}$, we have : 
\begin{equation*}
\begin{split}
   \phi(q) = 
  \min\biggl( &  \min_{i,j,k} \Bigl\lbrace \sqrt{RT(i,j,k)} \text{ such that } 
\|  \vit(q)-\um(i,j,k)   \|\leq c\sqrt{RT(i,j,k)}\Bigr\rbrace \biggr. , \\
& \, \biggl. \max_{i,j,k} \sqrt{RT(i,j,k)} \biggr)
\end{split}
\end{equation*}
\end{prop}
This means that if we take a velocity $\vit(q)$ of the fine grid, then
among all the distribution functions whose support contains $\vit(q)$,
one of them has a smallest support, and $2c\phi(q)$ is the size of
its support. This support function hence tells us how the fine Cartesian
grid should be refined, or coarsened, around $\vit(q)$.

The generation of the corresponding adapted grid is described in the
following section.

\begin{remark}
It may happen that with such a strategy, the support function does not
take sufficiently into account the support of the wall Maxwellian (for
diffuse reflection conditions on a solid obstacle): this Maxwellian is
centered around $\um=0$ and has a temperature $T_{wall}$ which is often
the smallest temperature of the computational domain. This can be
the case when the macroscopic CNS flow is not resolved enough close to
the wall, or if the CNS equations are solved with slip boundary
conditions. In this case, it is safer to add the element
$\phi(q_{max}+1):=\sqrt{RT_{wall}}$ in $\phi$ and to add the corresponding velocity
$\vit(q_{max}+1):=0$ in the fine grid.
\end{remark}

   \subsection{AMR grid generation}
   \label{subsec:algo_amr}

The idea is to start with a unique cell defined by the bounds of the
fine velocity grid (that is the full velocity domain), and then to apply a recursive algorithm: each cell
is cut if its dimensions are larger than the minimum of $a\phi$ in the
cell. The algorithm is the following:
\begin{algo}[Recursive refinement of a cell ${\cal C}$]
\label{algo:AMR}
$\ $

\begin{enumerate}
  \item compute the minimum of the field $a\phi$ in the cell ${\cal
      C}$, that is to say the minimum of the elements of $a\phi$ that have the
    same indices as the discrete velocities of the fine grid 
    included in this cell ${\cal C}$:
\begin{equation*}
  m:=a\min\lbrace \phi(q), \text{ such that } \vit(q)\in \text{cell } {\cal C}\rbrace
\end{equation*}  
\item if one edge of ${\cal C}$ is larger than $m$, then the cell is
  cut into 8 new subcells on which the refinement algorithm is applied, recursively.
  \item else, the cell is kept as it is, and the cell and its vertices
    are numbered.
\end{enumerate}
\end{algo}

At the end of the algorithm, the final grid satisfies the following
property: every cell has a size smaller than $a\times$the
minimum of the support function in the cell.

   \subsection{An example}

We anticipate on the numerical results that will be discussed in
section~\ref{sec:res} to illustrate the previous algorithms. 

\begin{figure}[h]
    \centering
    \includegraphics[width=0.5\textwidth]{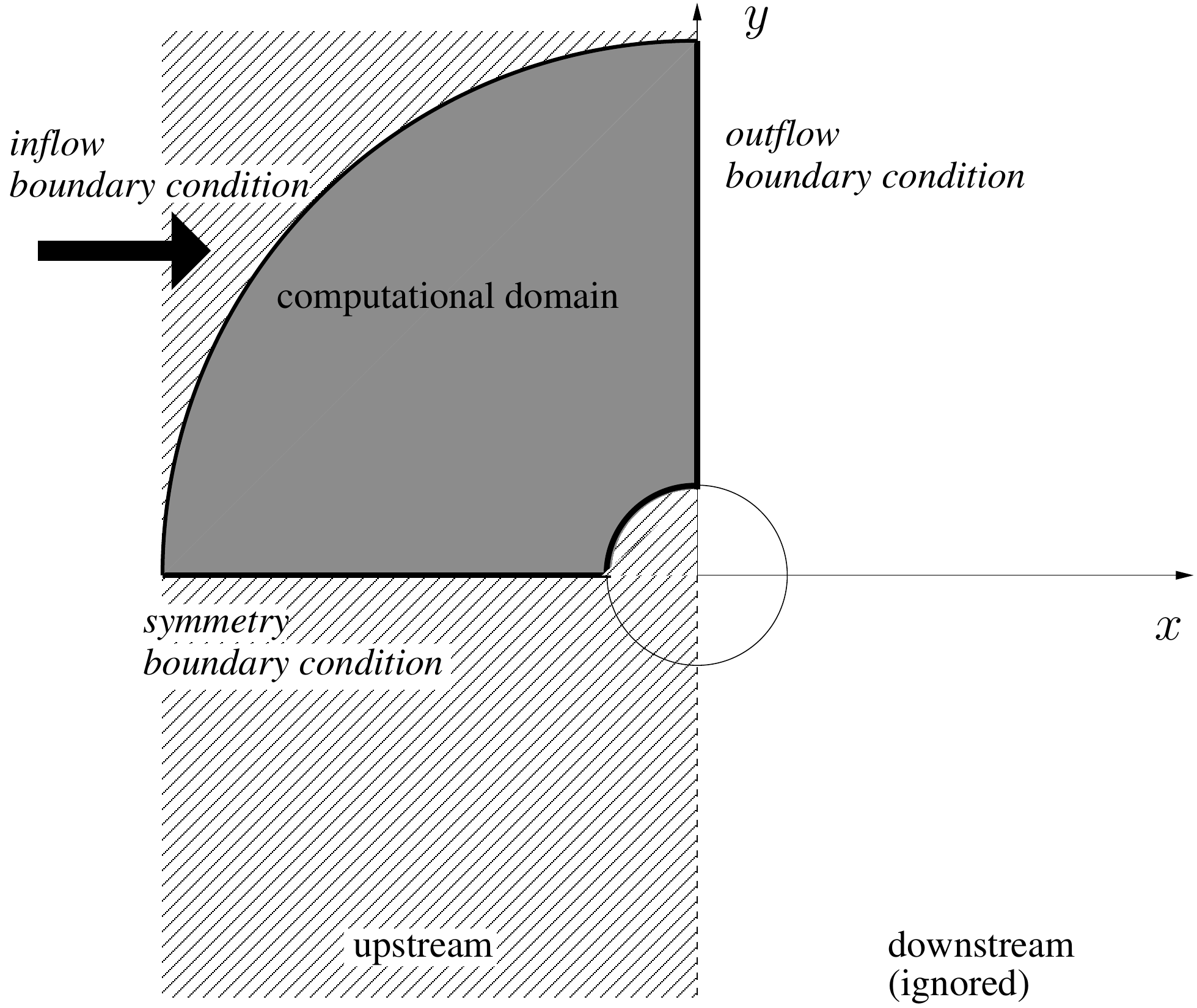}
    \caption{\modif{Plane flow around a cylinder: geometry and
      computational domain. By symmetry
      with respect to the axis $x$, the computational domain is
      defined for $y>0$ only. The downstream flow is not simulated.}}
    \label{fig:geom_cylindre}
  \end{figure}

\begin{figure}
  \centering
    \includegraphics[height=0.4\textheight]{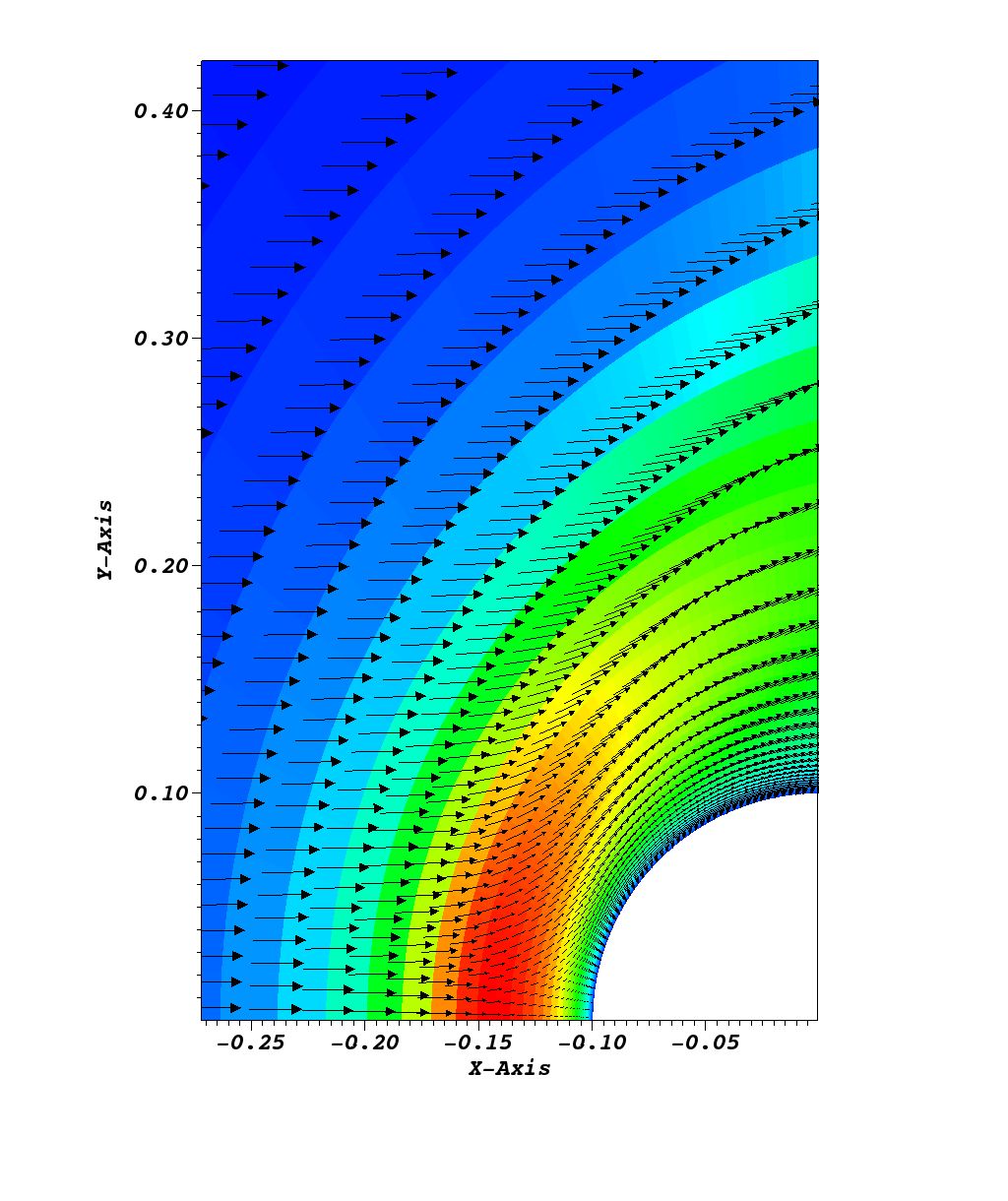} 
\includegraphics[height=0.4\textheight]{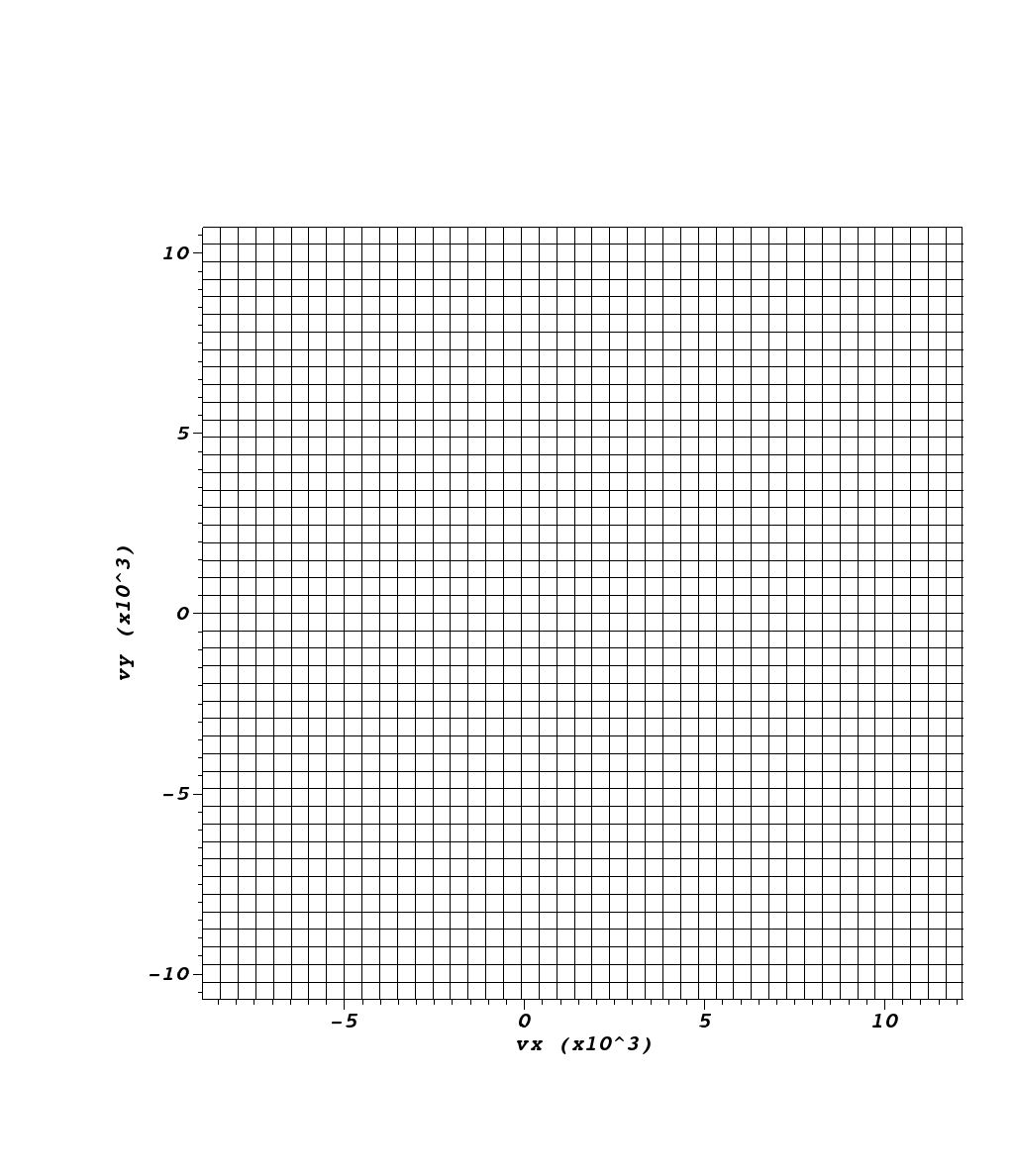}
    \caption{\modif{Plane flow around a cylinder:} CNS velocity and temperature fields (left), corresponding
      fine Cartesian
      velocity grid (right). }
    \label{fig:CNS_grid}
 \end{figure}

\begin{figure}
   \centering
    \includegraphics[width=0.45\textwidth]{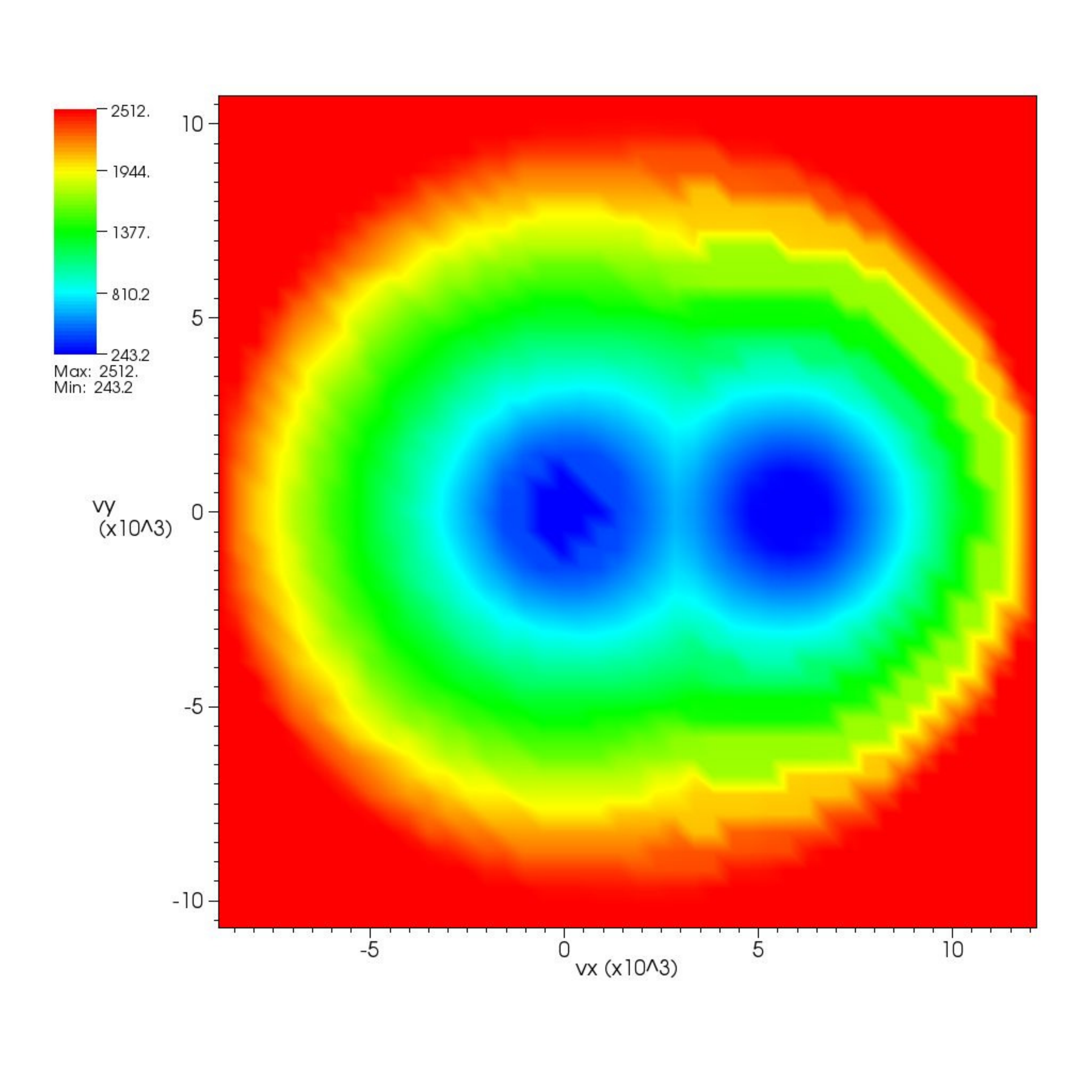} \includegraphics[width=0.45\textwidth]{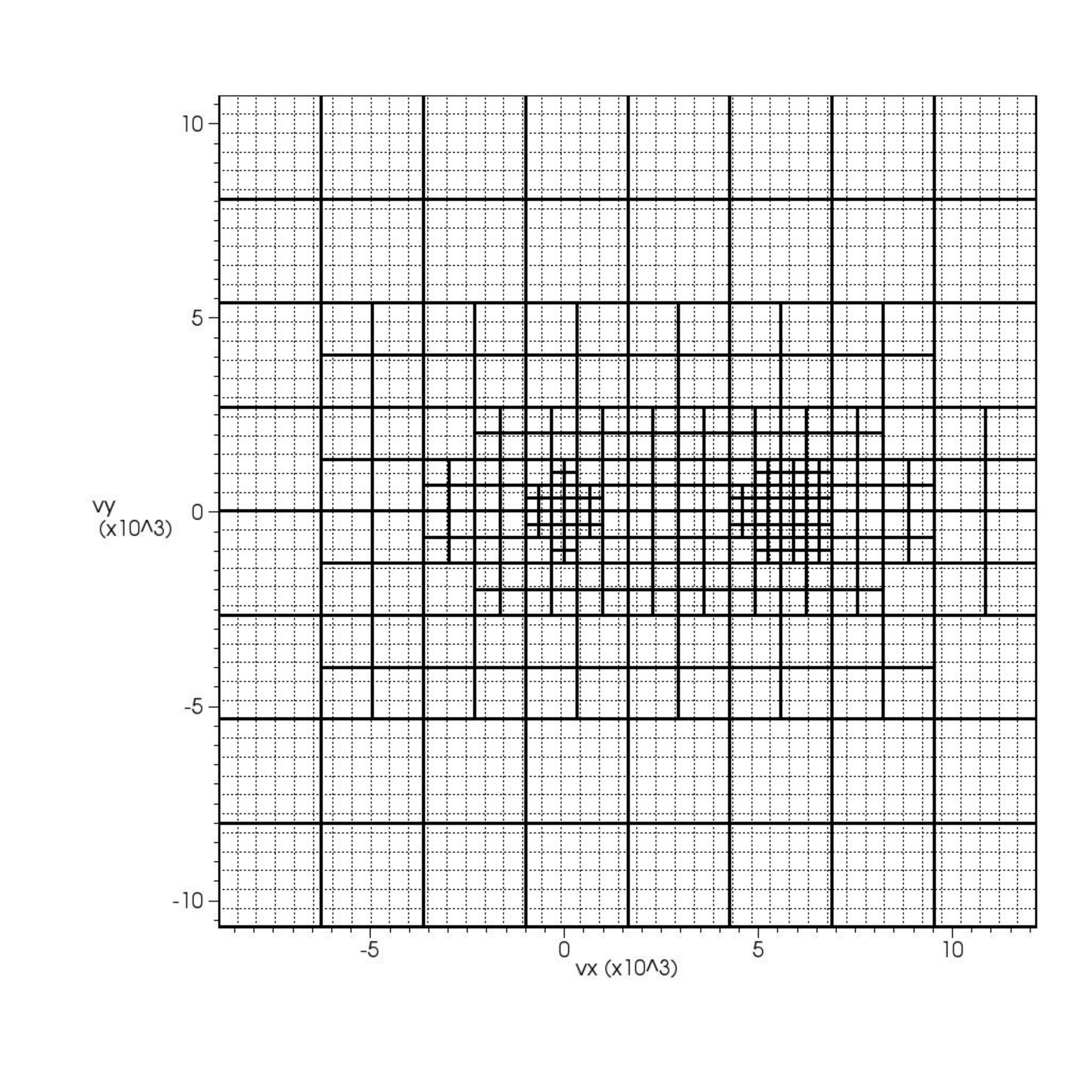}
   \caption{\modif{Plane flow around a cylinder.} Left: support function, Right: velocity grids (solid line:
     induced AMR velocity grid, dotted line: initial fine Cartesian
     grid).}
   \label{fig:AMR_grid}
 \end{figure}

The test case is a \modif{2D} steady flow around a cylinder at Mach 20 \modif{(see the
geometry in figure~\ref{fig:geom_cylindre} and the physical parameters
in section~\ref{subsec:2D})}. A CNS
pre-simulation and the use of
formula~(\ref{eq-bounds_CNS})-(\ref{eq-step_CNS}) with parameters $c=4$
and $a=2$ give us a fine velocity grid with $52\times 41$ points, see
figure~\ref{fig:CNS_grid}.  Algorithms~\ref{algo:support}
and~\ref{algo:AMR} applied to the CNS fields give the support function
and the AMR velocity grid shown in figure~\ref{fig:AMR_grid}. The AMR
velocity grid has 529 points, which is one fourth
as small as the original fine Cartesian grid. 

Note that the AMR grid is refined around the zero velocity and the
upstream velocity: these velocities correspond to the flow close to
the solid boundary and to the upstream boundary, where the temperature
is low, thus where the distribution functions are narrow, and hence
where, indeed, we need a fine velocity resolution. At the contrary, the grid is
coarse close to its boundaries: the discrete velocities
are very large there and are all in the tails of all the distribution
functions, where, indeed, we
do not need a fine resolution.

The accuracy of the computation on this AMR grid is studied in section~\ref{sec:res}.

   \subsection{Discrete Velocity Approximation on the AMR velocity grid}
When one wants to transform a standard discrete velocity method based
on a Cartesian grid to a method using our AMR grid, two points have to
be treated carefully: the computation of the moments of $f$, and the
approximation of the collision operator. In this paper, we only use the
BGK collision model of the Boltzmann equation: therefore its
approximation (based on the conservative approach of~\cite{luc_jcp})
reduces to the problem of a correct approximation of the moments of
$f$, as it is described below.

First note that the standard discrete velocity approach consists in
replacing the kinetic equation~(\ref{eq-kin}) by the following
discrete kinetic equation
\begin{equation}
\dt f_q + \vit_q \cdot \nabla_\x f_q = Q_q(f),
\end{equation}
where $f_q(t,x)$ is supposed to be an approximation of $f(t,x,\vit_q)$
for every discrete velocity $\vit_q$, with $q=0$ to $q_{max}$. The
discrete distribution function is $f=(f_q(t,x))$, and $Q_q$ is
the discrete collision operator. The moments of the discrete
distribution function are obtained by using a quadrature formula to
replace~(\ref{eq-mts}) by
\begin{equation*} 
(\rho(t,\x),\rho\um(t,\x),E(t,\x))
=\sum_{q=0}^{q_{max}}(1,\vit_q,\demi|\vit_q|^2)f_q(t,\x)\omega_q,
\end{equation*}
where $(\omega_q)$ are the weights of the quadrature.

For the BGK model~(\ref{eq-bgk}), the conservative
approximation of~\cite{luc_jcp} gives the following discrete BGK
collision operator
\begin{equation*}
  Q_q(f)=\frac{1}{\tau}(\Maxw_q[\rho,\um,T]-f_q),
\end{equation*}
where $\Maxw_q[\rho,\um,T]$ is the discrete
Maxwellian that has the same discrete moments as $f$. It can be shown that there exists $\alpha\in\R^5$
 such that $\Maxw_q[\rho,\um,T]=\exp(\alpha\cdot m(\vit_q))$, that is to say,
$\alpha$ is the unique solution of the nonlinear system
\begin{equation*}
  \sum_{q=0}^{q_{max}} m(\vit_q)\exp(\alpha\cdot m(\vit_q)) \omega_q = 
  \begin{pmatrix}
    \rho \\ \rho \um \\ E
  \end{pmatrix},
\end{equation*}
where we note $m(\vit)=(1,\vit,\demi |\vit|^2)^T$. This system can be
solved by a Newton algorithm (see details in~\cite{luc_jcp} and a more
economic version of the algorithm due to Titarev in~\cite{T_cicp}).

Consequently, to apply this approach to our AMR velocity grid, we just
have to chose a quadrature formula. First, we propose the standard
$\mathbb{Q}_1$ bilinear interpolation of $f$ on each cell of the AMR grid (by
using the four corners of the cell). In this case, the quadrature
points are the nodes of the grid, and so are the discrete velocities
that are used in the transport term. The quadrature weights are
\begin{equation*}
\omega_q=\frac{1}{2^d}\sum_{\text{cell } {\cal C}_l \ni \vit_{q}}|{\cal C}_l |,  
\end{equation*}
where $d=3$ in 3D or $2$ in 2D, and $|{\cal C}_l |$ is the volume (or
the area in 2D) of ${\cal C}_l $.

The number of discrete velocities can be decreased if we take a
simpler $\mathbb{P}_0$ (constant per cell) interpolation. Here, the
quadrature points are the centers of the cells, and so are the
discrete velocities $\vit_q$ that are used in the transport term. The weights
are simply the volume (or the area) of the cells: $\omega_q=|{\cal
  C}_q |$. Note that the number of discrete velocities that are used
with this approach is equal to the number of cells of the AMR grid,
which is smaller than the number of nodes, see
section~\ref{subsec:2D}. We advocate the use of this latter approach, since it
allows to save a lot of computer memory, and since we observed with our tests that the
$\mathbb{Q}_1$ interpolation does not give more accurate results.

\begin{remark}
  Our AMR velocity grids are often very coarse at their boundaries
  (since velocities are very large there, and are in the tails of all
  the distribution functions). Consequently, passing from
  $\mathbb{Q}_1$ to a $\mathbb{P}_0$ grid may lead to a grid which is
  not large enough. Indeed, if the $\mathbb{Q}_1$ grid is of length
  $L$ (in the $x$ direction, for instance) and its outer cells are of
  length $l$, then the length of the $\mathbb{P}_0$ grid now is
  $L-l$. If $l$ is large, the $\mathbb{P}_0$ grid is not
  large enough. In that case, it is thus safer to modify step 2 of
  Algorithm~\ref{algo:support} by increasing the bounds of the fine
  Cartesian grid: $v^{\alpha}_{max}$ and $v^{\alpha}_{min}$ are
  replaced by $v^{\alpha}_{max}+\Dv$ and
  $v^{\alpha}_{min}-\Dv$ for example, where $\Dv$ is the step of the fine grid.
\end{remark}

   \subsection{AMR grid generation for axisymmetric flows}
   \label{subsec:axi}

Many interesting 3D flows have symmetry axis (like flows around
axisymmetric bodies with no incidence). In that case, if we assume
that the symmetry axis is aligned with the $x$ coordinate, it is
interesting to write the kinetic equation in the cylindrical
coordinate system $(x,r,\varphi)$, since the distribution function $f$
is independent of $\varphi$, and we get
\begin{equation} \label{eq-bgkaxi}
  \dt f + v_x\dx f + \zeta \cos \omega \partial_r f -
  \frac{\zeta\sin\omega}{r} \partial_\omega f = Q(f),
\end{equation}
where $\zeta$ and $\omega$ are the cylindrical coordinates
of the velocity in the
local frame $(e_r,e_\varphi)$, that is to say $v\cdot
e_r=\zeta\cos\omega$ and $v\cdot e_\varphi=\zeta\sin\omega$. 

If the velocities of the upstream flow and of the solid boundaries have no component in the $\varphi$
direction, then $f$ is even with respect to $\omega$. Indeed, using the
Galilean invariance of equation~(\ref{eq-bgkaxi}), we apply the
symmetry with respect to the $(x,z)$ plane to the equation to get that
$f(t,x,r,v_x,\zeta,-\omega)$ is also a solution. 

In this case, the discrete velocity grid must be constructed for
variables $(v_x,\zeta,\omega)$, where $\zeta$ is non negative,
$\omega$ is in $[0,\pi]$, and $v_x$ can take any value. The
construction of this grid is slightly different from the Cartesian
case. We briefly comment on these differences in this section.

First, bounds in $v_x$ and $\zeta$ directions can be easily obtained
and we set
\begin{equation}  \label{eq-bounds_steps_axi}
\begin{aligned}
 v_{x,max}& =\max_{\Omega}\lbrace u_x^{CNS}+c\sqrt{RT} \rbrace, 
&v_{x,min}& =\min_{\Omega}\lbrace u_x^{CNS}-c\sqrt{RT} \rbrace ,   \\
 \zeta_{max}& = \max_{\Omega}\lbrace u_r^{CNS}+c\sqrt{RT} \rbrace, 
& \zeta_{min}& = 0,
\end{aligned}
\end{equation}
where $u_x$ and $u_r$ are the macroscopic velocities in the axial and
radial directions. Since $\omega$ is a bounded variable, there is
nothing else to do at this level.

The steps $\Delta v_x$ and $\Delta \zeta$ of the initial grid are
obtained by~(\ref{eq-step_CNS}) for $v_x$ and $\zeta$. However, the
choice of the step in the $\omega$ direction is not obvious. Indeed,
this variable has no link with the width of the
distribution. Moreover, the grid in this direction must be fine enough
to capture the derivative $\partial_\omega f$, while in the other
directions, only some moments have to be captured. After several tests
(and from the experience made in~\cite{luc_jcp}), it seems that a grid
step $\Delta \omega$ such that there are 30 points uniformly distributed in $[0,\pi]$ is a
good compromise between cost and accuracy.

To refine/coarsen this initial cylindrical grid, we choose to use the
procedure described in the steps of sections~\ref{subsec:algo_support}
and~\ref{subsec:algo_amr} in the plane $(v_x,\zeta)$ to generate a
two-dimensional AMR grid. Indeed, for the same reasons as mentioned
above, there is no reason to refine or coarsen the grid in the $\omega$
direction. Then, the complete grid is obtained by a rotation of the
two-dimensional AMR grid in the $\omega$ direction, with the same step
$\Delta \omega$ as in the initial grid. This procedure is simple and
turned out to be very efficient.

Since the AMR grid is uniform in the $\omega$ direction, the
discretization of~(\ref{eq-bgkaxi}) is made with the approximations
derived in~\cite{luc_jcp}.

\section{Numerical results}
\label{sec:res}

   \subsection{The 3D kinetic code}

   Our 3D kinetic code is written in Fortran, and can deal with 2D
   planar flows, or 2D axi-symmetric flows as well as 3D flows on
   multiblock curvilinear structured grids. A steady state solution of the BGK
   equation for monoatomic or polyatomic gases is computed with a
   linearized implicit finite volume scheme, based on a velocity
   discretization which is conservative and entropic (see
   \cite{luc_m3as}). In axisymmetric cases, a specific scheme (named
   T-UCE) is used in order to ensure a conservative and entropic
   discretization of the transport operator (see \cite{luc_jcp}).

In order to deal with high computational costs (in 3D configurations for
example), an hybrid parallel implementation is used: a domain
decomposition in the position space used with the library routines
specified by the Message Passing Interface (MPI) can deal with a large
number of mesh cells. Note that this is the opposite of the strategy
used in~\cite{T_cicp} (which uses a domain decomposition in the
velocity space). Moreover, the OpenMP library is used with 8 threads
for loops computations that are local w.r.t the position (computation
of the moments, Maxwellian, etc.), thus large numbers of velocity
points can be reached. We also use OpenMP for loops that are local
w.r.t the velocity, like the computation of the transport terms, for
instance. The communication between different domains are treated as
explicit boundary conditions. This code has been run on the
super-computer TERA-100 of the CEA (that has around $140 \, 000$ cores Intel Xeon
7500), by using at most 1600 cores. 

Since a general description of our code is not the main subject of
this paper, we do not give more details here. However, for
completeness, note that interested readers can find a more detailed
description of the code in the appendix.

The distribution is initialized by using a CNS
pre-simulation (also used for the AMR strategy): the number of
iterations before reaching convergence is then considerably
decreased (see section~\ref{subsec:init}).

For the 2D plane examples shown below, we use a simpler version of the code,
called CORBIS. This code is based on the same tools but makes use
of the reduced distribution technique to reduce the velocity space to
2D only. This code only uses OpenMP directives to work on parallel
computers (see~\cite{ADM2009}). It has been used on a SMP node of 48 cores
AMD-Opteron-8439-SE. 

\modif{\label{page:code_cns}For the compressible Navier-Stokes simulations, we use a 3D code
developed by the CEA-CESTA for hypersonic flows at moderate
to low altitudes. The same techniques as in our 3D kinetic code are
used: a finite volume linearized implicit scheme with multiblock
structured meshes. On solid walls, no-slip and isothermal boundary
condition are used ($u=0$ and $T=T_{wall}$).}

   \subsection{A 2D plane example}
   \label{subsec:2D}

   Here, we illustrate our approach with the simulation of a steady
   flow over a infinite cylinder of radius $0.1$m at Mach 20, for
   density and pressure of the air at an altitude of 90 km. The gas
   considered here is argon (molecular mass$=6.663\times 10^{-24}
   kg$). Namely, we have $\rho=3.17\times 10^{-6} kg/m^3$,
   $u=5.81\times 10^3 m/s$, and $T=242.4 K$. \modif{The temperature of
     the wall is $293 K$.}  \modif{The downstream flow is ignored, and
     outflow boundary conditions are used at the boundary of the right
     side of the domain. Since the flow is invariant along $z$, the
     simulation is made with the 2D plane model in the plane
     $(x,y,0)$. Finally, since the flow is symmetric with respect to
     the line $y=0$, the computation is made on the upper part on this
     plane only (see figure~\ref{fig:geom_cylindre}). See appendix~\ref{subsec:num_bc} for
     the implementation of the boundary conditions in the code.}

The space mesh uses $50\times 50$ cells, with a refinement such that
the first cell at the solid boundary is smaller than one fifth of the mean
free path at the boundary.

A CNS pre-simulation and the use of
formula~(\ref{eq-bounds_CNS}-\ref{eq-step_CNS}) with parameters $c=4$
and $a=2$ give us a fine velocity grid with $52\times 41$ points
(bounds $\pm 11\,500$ and $\pm 8\,900$ for $v_x$ and $v_y$ with a step
of $450$, in m.s$^{-1}$), see figure~\ref{fig:CNS_grid}.
Algorithms~\ref{algo:support} and~\ref{algo:AMR} applied to the CNS
fields give the support function and the AMR velocity grid shown in
figure~\ref{fig:AMR_grid}.  The AMR velocity grid has 529 points,
which is one fourth as small as the original fine Cartesian grid. This
gain can be further increased if we define the discrete velocities as
the centers of the cells of the AMR grid rather than its
vertices. This gives 316 discrete velocities, hence with a gain of 6.7
times instead of 4.

First, we compare the CPU time required by the code CORBIS with the different velocity grids. While the number of iterations
to reach steady state is approximately the same with both grids,
the CPU time required by the original fine Cartesian velocity grid is around 7 times as large as with the
new AMR grid, which is almost the same ratio as the ratio of the number of
discrete velocities. 
The memory required with the uniform grid method is around 170 MB whereas with the use of the AMR grid, only 25 MB 
of memory storage is required.
This shows that the new AMR grid leads to a real gain
both in memory storage and CPU time.

\begin{figure}
  \centering
  \includegraphics[width=0.65\linewidth,clip=true]{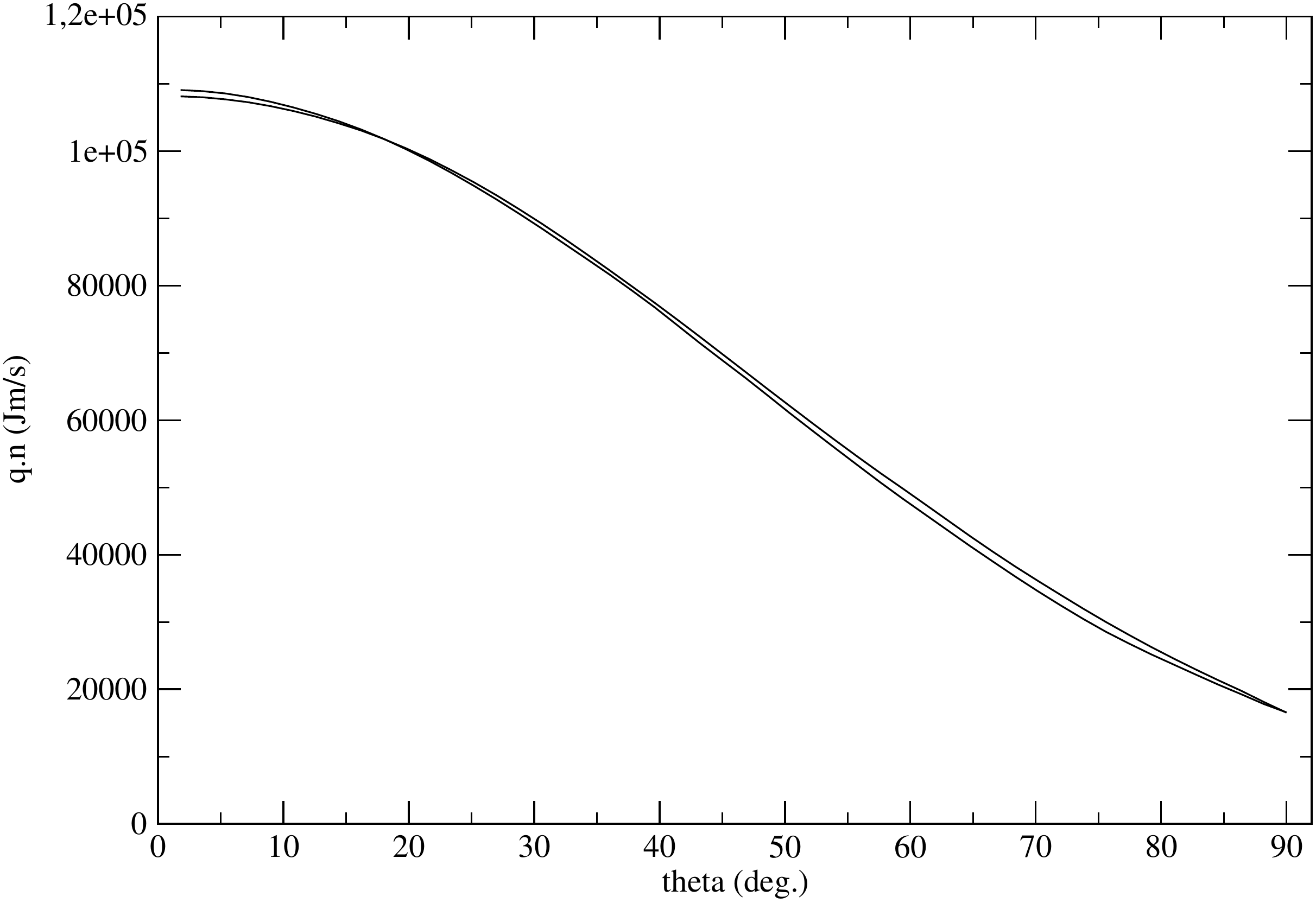} 

\vspace{10ex}

  \includegraphics[width=0.6\linewidth,clip=true]{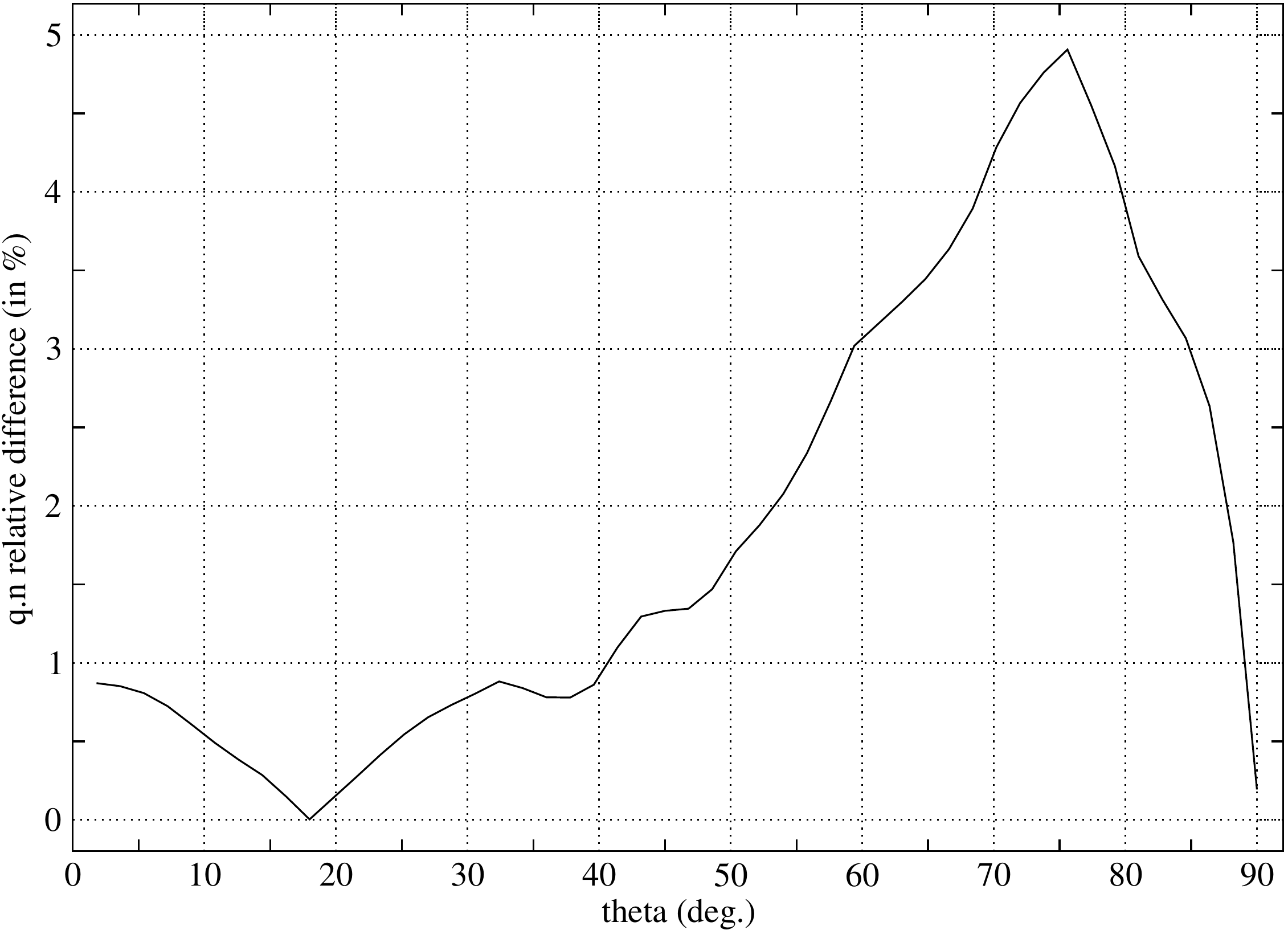}
  \caption{\modif{Plane flow around a cylinder.} Comparison of the
    component of the heat flux normal to the
    solid wall: fluxes obtained with the fine Cartesian grid and the
    AMR grid (top), relative difference in percent
    (bottom). \modif{The angle theta defines the position of the
      point along the solid boundary, such that theta$=$0 at the
      stagnation point and theta$=$90 at the end of the wall}.}
  \label{fig:comp_flux_amr}
\end{figure}

Then, we compare the accuracy of the results with the two grids for
the macroscopic quantities. We compute the normal component of the heat flux to the boundary,
which is a quantity of paramount importance in aerodynamic
simulations: we find a maximum relative difference lower than 5\%,
which is reasonably small (see the profile of this flux on
figure~\ref{fig:comp_flux_amr}). We also compute the differences for
the density, temperature and pressure in the whole computational domain:
\begin{itemize}  

\item the mean quadratic relative difference over all the cells of the
    computational domain is 5\% for the density, 1\% for the
    temperature, 0.6\% for the horizontal velocity, and 1.2\% for the
    vertical velocity;

  \item the maximum relative difference on each cell of the
    computational domain is 40\% for the density, 69\% for the
    temperature, and 157\% for the velocity.
\end{itemize}
The maximum relative differences are observed at the solid wall for
the density and the velocity, and in the upstream flow for the
temperature. 

This difference is quite large, and can be explained as follows. The
smallest cells of our AMR velocity grid (that are around small
velocities, like velocities at the solid wall) turn out to be smaller than
the cells of the Cartesian grid (size $330$ instead of $450$). This
means that our results with the AMR grid are probably {\it more
  accurate} than the results of the Cartesian grid. Consequently, the
Cartesian grid results should not be considered as a reference for this
comparison. 

To confirm this analysis, we make a new simulation with a
Cartesian grid with a uniform step of $330$ (like the smallest step of
the AMR grid). We observe that relative differences between the
Cartesian grid and the AMR grid are much smaller: 
\begin{itemize}
  \item the mean quadratic relative difference over all the cells of the
    computational domain is 2\% for the density, 0.5\% for the
    temperature, 0.3\% for the horizontal velocity, and 0.7\% for the
    vertical velocity; 

  \item the maximum relative difference on each cell of the
    computational domain is 12\% for the density, 13\% for the
    temperature, 80\% for the velocity.
\end{itemize}
The maximum relative difference is still too large for the velocity
(80\% at the boundary), but this quantity is very small in this zone and
probably requires smaller velocity cells (that is a smaller parameter
$a$ for both grids). However, this inaccuracy does not deteriorates
the results on the other quantities, in particular for the heat flux
at the boundary. Indeed, the comparison of the heat fluxes is even
better, since we find a relative difference lower than 2.5\%, which is
excellent. The comparison in terms of CPU time and memory storage is
of course more favorable to the AMR grid here.

Note that we also did the same computations with an finer AMR grid
obtained with the parameter $a=1$: the difference with $a=2$ is less
than 2\% for the heat flux at the solid wall, which is quite
good. The value $a=2$ is then clearly sufficient. However, if one is interested with the
flow velocity at the boundary, we advocate to use $a=1$, since the
differences here are around 15\%.

   \subsection{A 2D axisymmetric example}
  \label{subsec:res_axi}

\begin{figure}[h]
    \centering
    \includegraphics[width=0.5\textwidth]{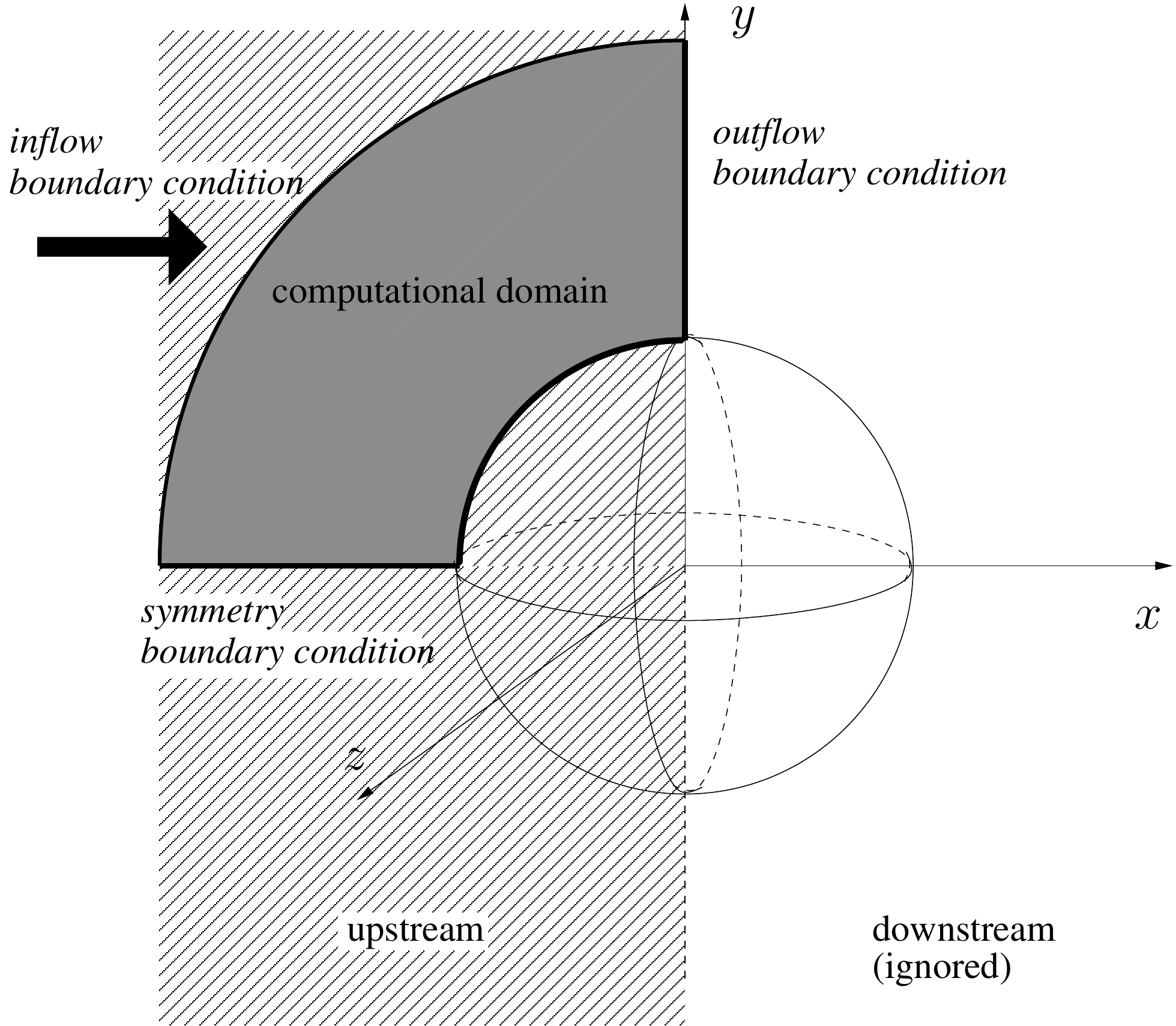}
    \caption{\modif{Axisymmetric flow around a sphere: geometry and
      computational domain. By rotational symmetry and by symmetry
      with respect to the axis $x$, the computational domain is
      included in the plane $(x,y,0)$. The downstream flow is not simulated.}}
    \label{fig:geom_sphere}
  \end{figure}

  Here we consider the flow around a sphere of radius 0.1 m, at Mach
  20, for density and pressure of the air at an altitude of 90 km
  \modif{(see figure~\ref{fig:geom_sphere})}.The gas is air (molecular
  mass$=4.81\times 10^{-26}$ kg), considered as a diatomic
  gas. Namely, the upstream density and pressures are $\rho=3.17\times
  10^{-6} kg/m^3$, and $p=0.16 N.m^{-2}$, and the temperature of the
  surface of the sphere is $T=280 K$. \modif{The downstream flow is
    ignored, and outflow boundary conditions are used at the boundary
    of the right side of the domain. Since the flow is rotationally symmetric,
    the simulation is made with the 2D axisymmetric model in the plane
    $(x,y,0)$. Finally, since the flow is symmetric with respect to
    the line $y=0$, the computation is made on the upper part on this
    plane only.}

The space mesh uses $70\times 50$ cells, with a refinement such that
the first cell at the solid boundary is smaller than one fifth of the mean
free path at the boundary. A CNS pre-simulation and the use of
formula~(\ref{eq-bounds_steps_axi}) with parameter $c=4$
 give us the following velocity bounds $-8\,700$ and $11\, 000$ for
 $v_x$, 0 and $7\,500$ for $\zeta$, and the bounds for $\omega$ are always
 0 and $\pi$, as explained in section~\ref{subsec:axi}.

All the following simulations were made with 140 domains and 4 OpenMP
threads per domain.

\begin{figure}[h]
    \centering
    \includegraphics[width=0.5\textwidth]{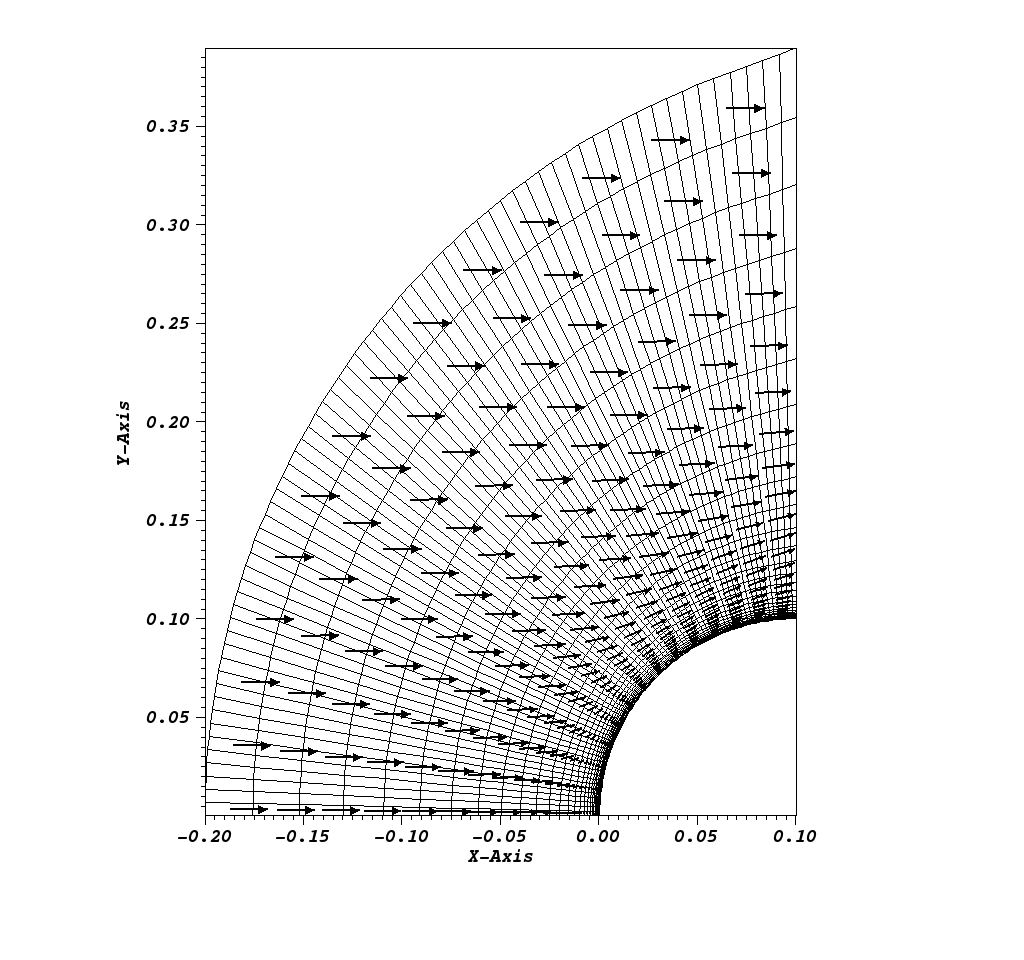}\includegraphics[width=0.5\textwidth]{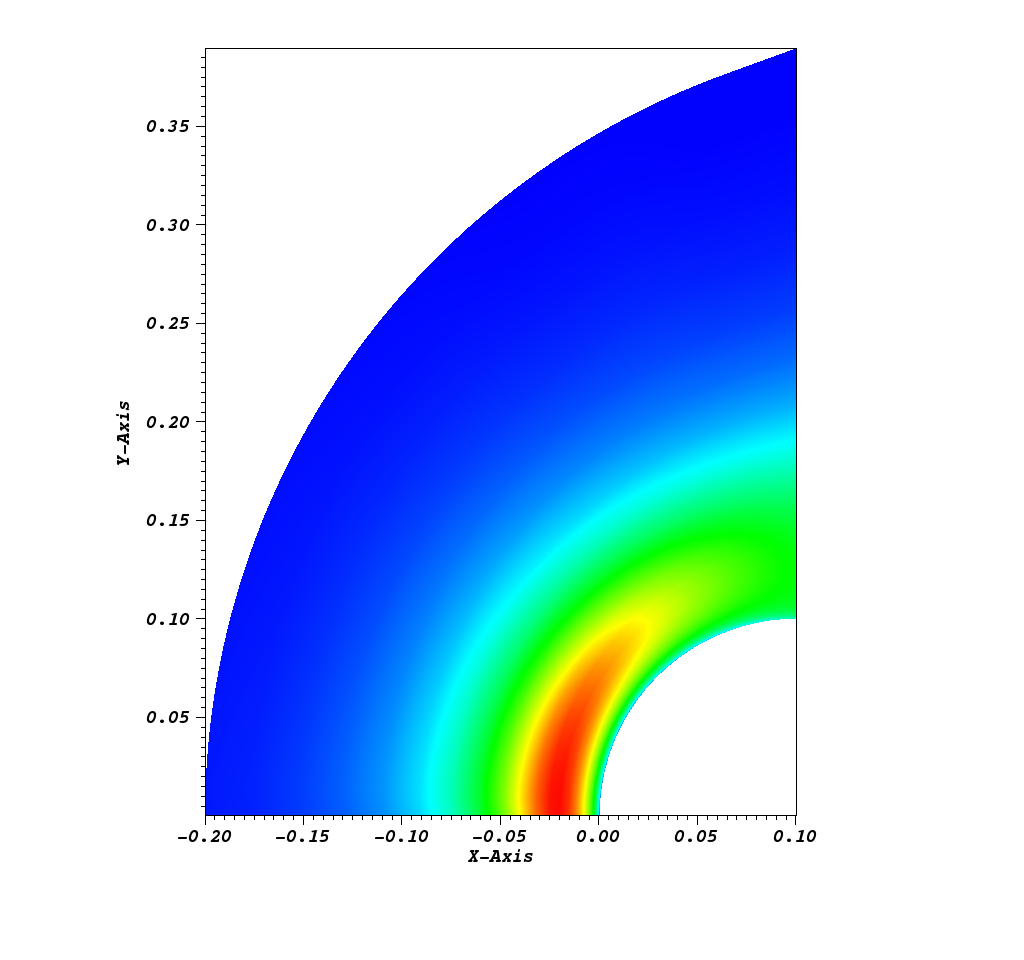}
    \caption{\modif{Axisymmetric flow around a sphere:} velocity field,
      \modif{spatial mesh} and temperature field.}
    \label{fig:fields_axi}
  \end{figure}

First, a fine uniform cylindrical velocity grid is obtained with
(\ref{eq-step_CNS}) and $a=2$: this gives $65\times33\times31=66\,
495$ discrete velocities. The steady state is reached after $1\, 828$ iterations,
in $5\, 476 \, s$. 
The corresponding temperature and velocity fields are shown in figure~\ref{fig:fields_axi}.

\begin{figure}[h]
    \centering
    \includegraphics[height=0.4\textheight]{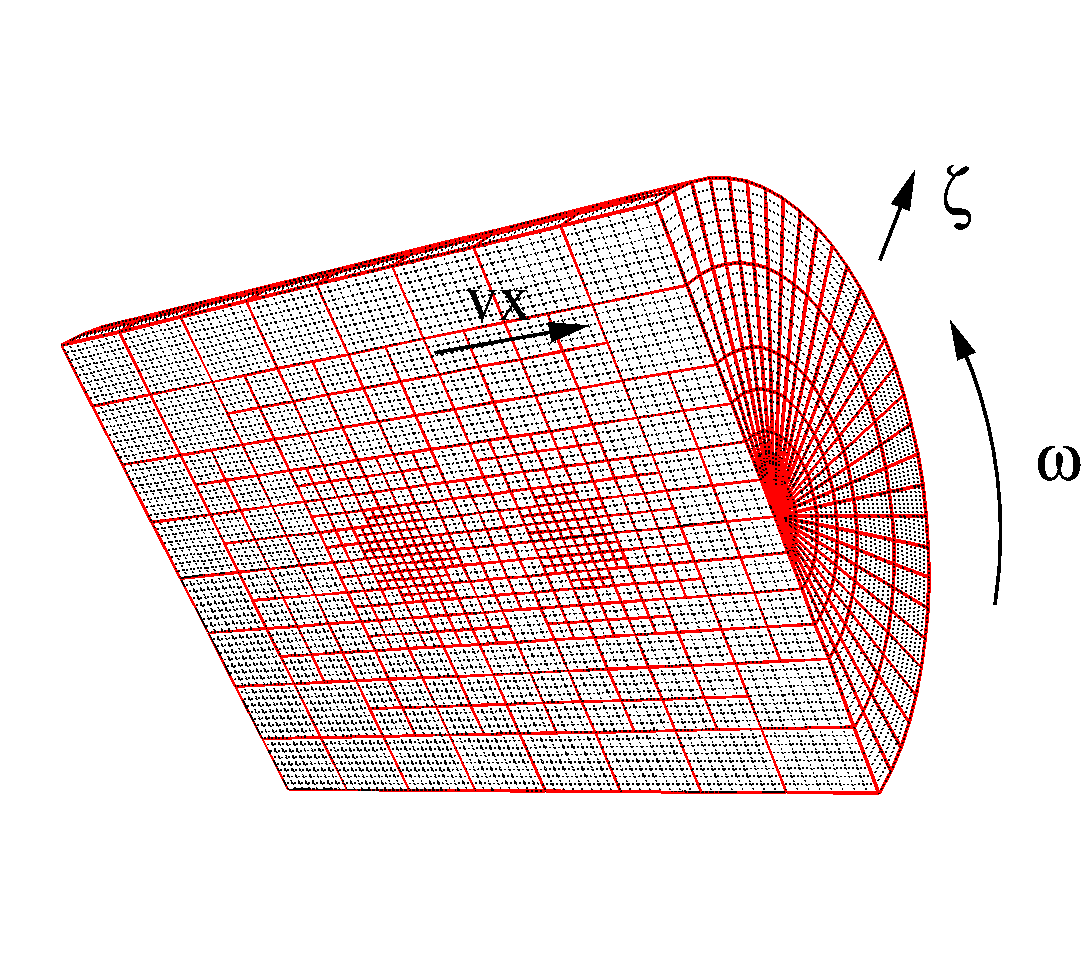}
    \includegraphics[height=0.4\textheight]{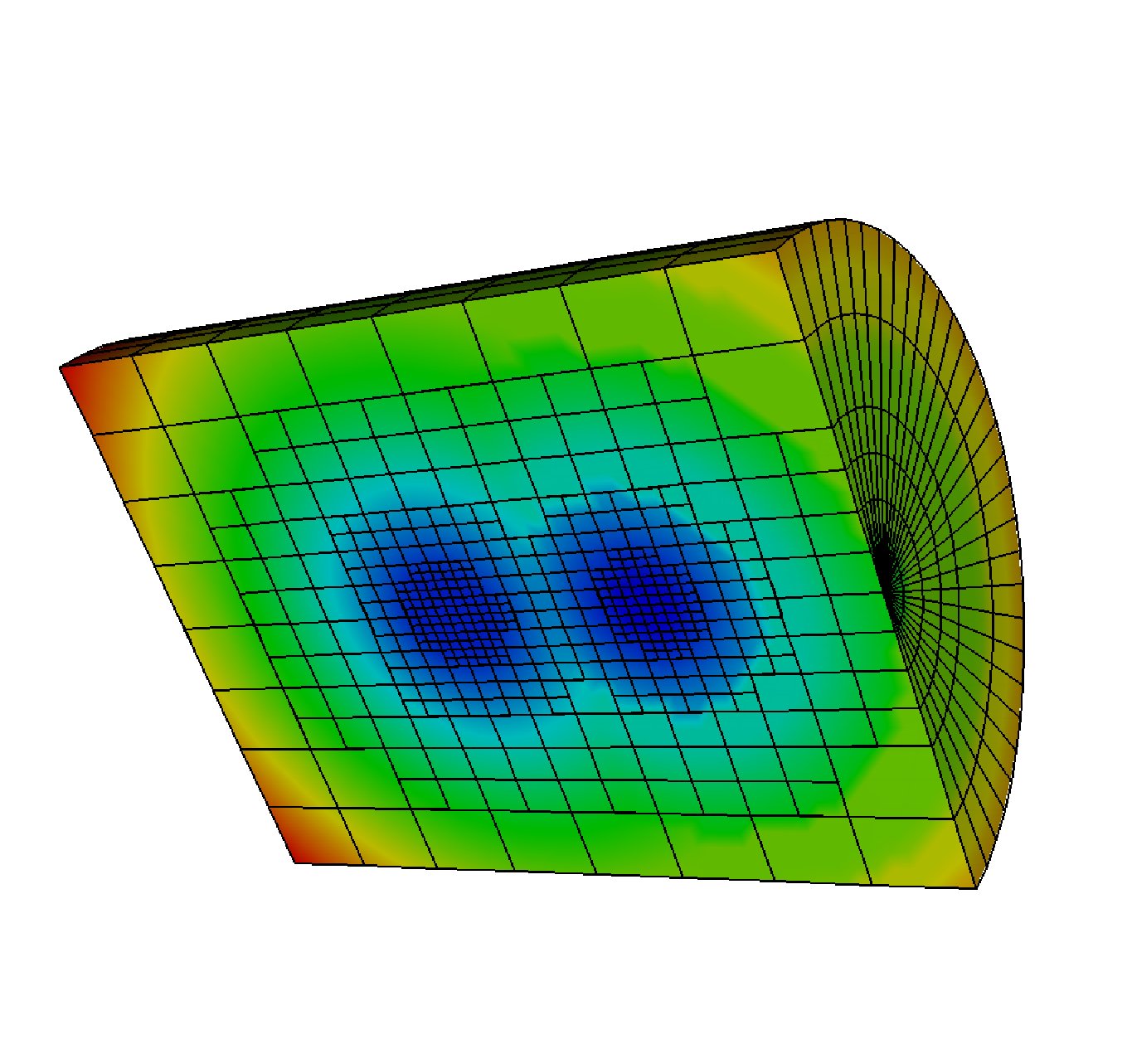}
    \caption{2D axisymmetric flow around a sphere. Top: uniform
      cylindrical \modif{velocity} grid (dashed lines) and AMR
      \modif{velocity} grid (solid) \modif{with the cylindrical velocity
      coordinate system $(v_x,\zeta,\omega)$}. Bottom:
      \modif{The same cylindrical AMR
      velocity grid with the contours of the support function based on
    the CNS pre-simulation}.}
    \label{fig:amr_cns_axi}
  \end{figure}

Then a cylindrical AMR grid is generated (with $a=2$): it gives $7\,
691$ discrete velocities (with the ${\mathbb Q}_1$ points), the steady
state is reached after $1\,899$ iterations, in $710 s$. The support
function and AMR grid are shown in figure~\ref{fig:amr_cns_axi}. Here,
the gain factors in memory and in CPU time with respect to the uniform
grid are around 9.

Now, we use the same AMR grid, but with the ${\mathbb P}_0$
points (that is to say the centers of the cells), which gives $6\,
900$ discrete velocities only.  The number of iterations is close ($1\, 796$), and the
CPU time is $562 s$. Here, the gain factors in
memory and CPU time, with respect to the uniform grid, are around 10.

\begin{figure}[h]
    \centering
     \includegraphics[width=\textwidth]{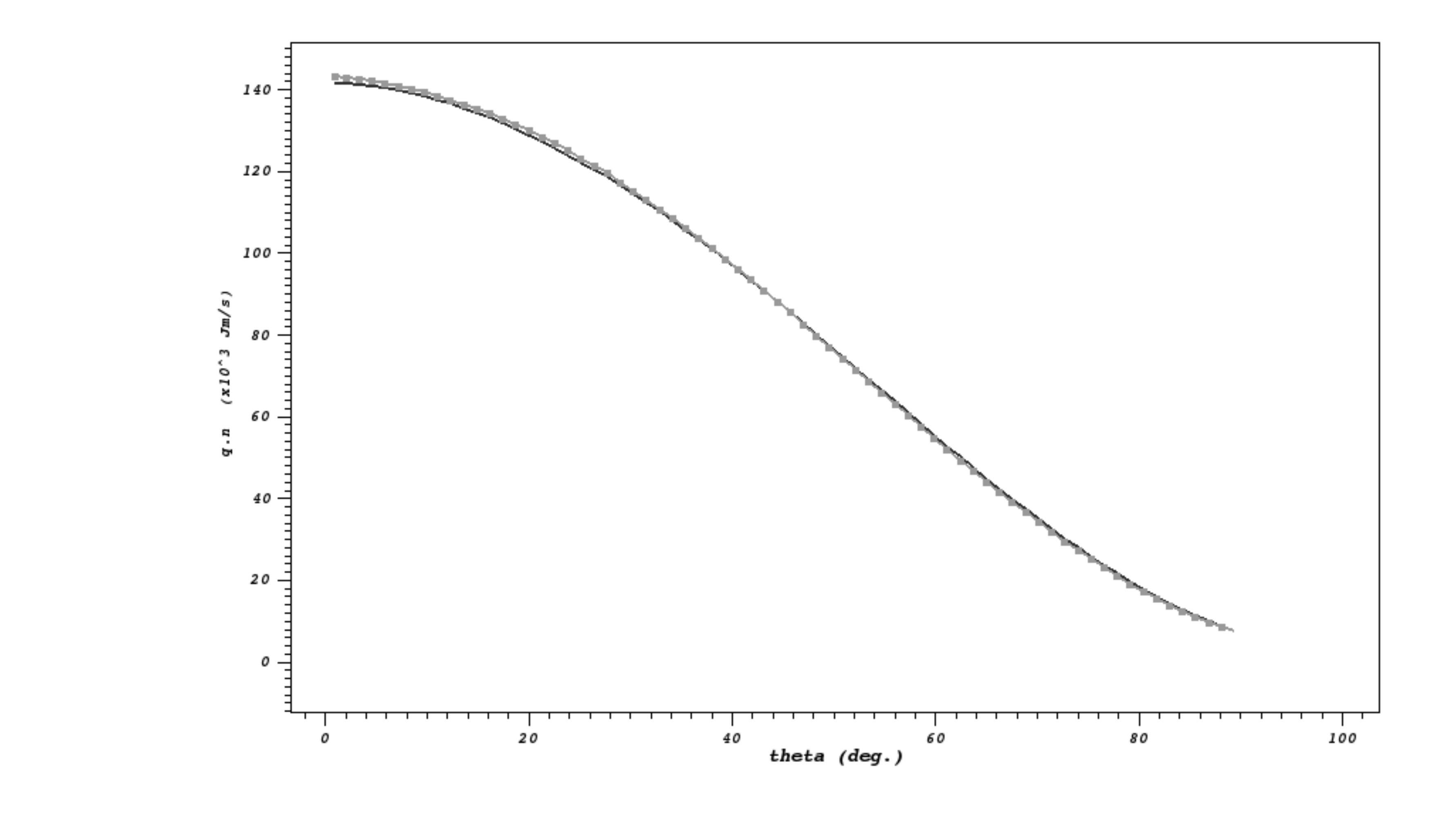}
    \caption{2D axisymmetric flow around a sphere: comparison of the
      component of the heat flux normal to
      the surface: fluxes obtained with the fine uniform
      \modif{velocity} grid (dot) and
      the AMR \modif{velocity} grid (solid)}.
    \label{fig:fluxes_axi}
  \end{figure}

Finally, we compare in figure~\ref{fig:fluxes_axi} the normal heat
flux at the boundary for the uniform grid and the AMR grid generated
by using the CNS fields (${\mathbb Q}_1$ version). As in the
2D example, we find that the fluxes are very close (with a relative
difference which is less than 1\%). We get the same results with the
${\mathbb P}_0$ version.

Consequently, the AMR grid (${\mathbb P}_0$ version) is the most efficient
here: it gives the same results as the fine uniform grid, for
computational and memory cost which is one ninth as small.

\begin{figure}[h]
    \centering
    \includegraphics[height=0.4\textheight]{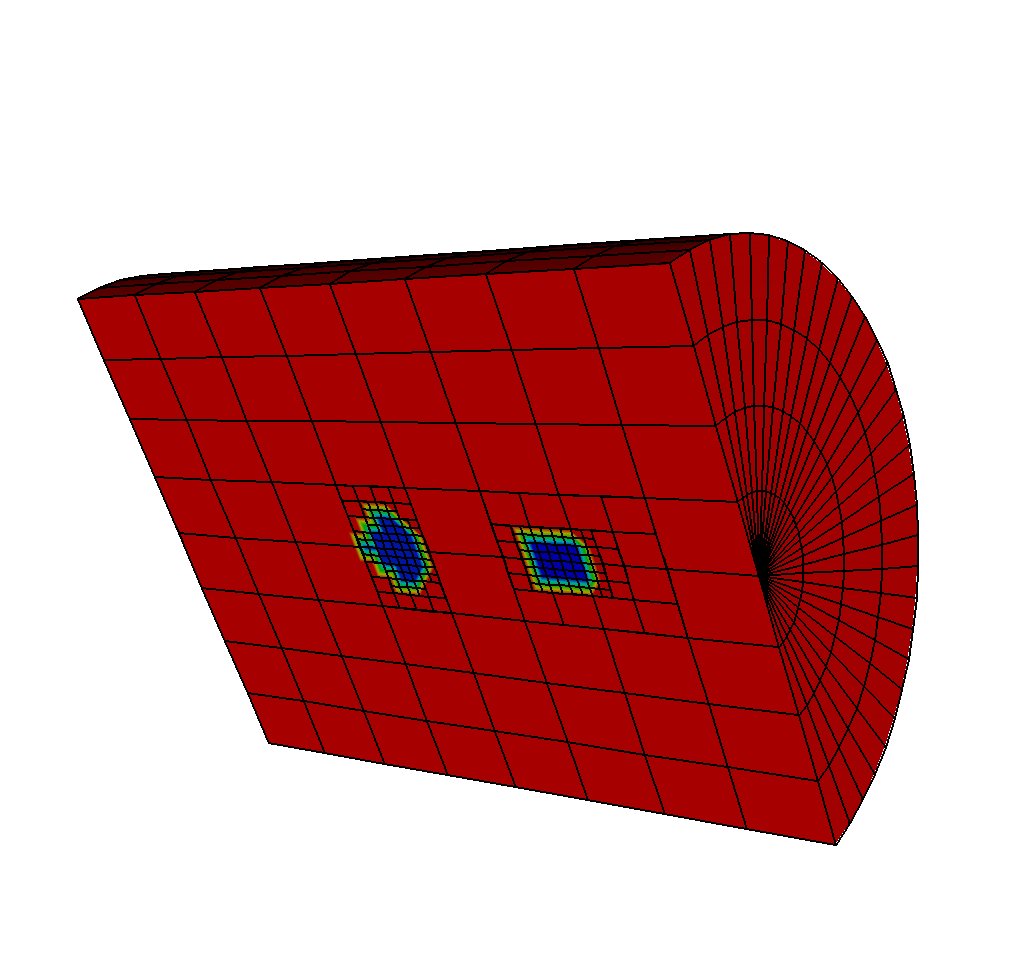}
   \caption{
     \modif{2D axisymmetric flow around a sphere: cylindrical AMR
       velocity grid computed with Rankine-Hugoniot relations
       rather that with the CNS pre-simulation, with contours of the
       corresponding support function. Note that the refined zones are
     much smaller than in the previous grid, and that the contours of
     the support function are different (see figure~\ref{fig:amr_cns_axi})}  }
    \label{fig:amr_rh_axi}
  \end{figure}

\begin{remark}
Before discussing a full 3D case, we briefly mention another
experiment we have done here. We have tried to see if the use of a CNS
pre-simulation could be avoided by using a direct estimation of the
extreme values of the velocity and the temperature given by the
Rankine-Hugoniot relations. Indeed, before the computation, we know
the upstream values $\um_{up}$ and $T_{up}$, and we can compute the
downstream values $(\um_{down},T_{down})=RH(\um_{up},T_{up})$ by the
Rankine-Hugoniot relations for a stationary shock. Then we have three
set of different values: upstream, downstream, and wall values
$(\um_{wall}=0,T_{wall})$, and we can use these three sets in our
algorithm~\ref{algo:support} instead of the CNS fields $(\um^{CNS},T^{CNS})$. 

With the same test case, this strategy gives a grid with  $3\, 962$ points
only (with $a=2$), which is half as less as with the CNS fields. The support function and AMR
grid are shown in figure~\ref{fig:amr_rh_axi}: of course, since
we only have three values for the macroscopic fields here, the grid is
not as smooth as with the CNS fields.

Moreover, the convergence to steady state is much slower (since we
cannot use the CNS solution as an initial state, see
section~\ref{subsec:init} for a comment on the initialization). More
important, the accuracy observed on the heat flux is not as good as
with the CNS fields: the difference  with the fine uniform grid is
here around 8\% instead of 1\%.

Consequently, even if strategy is simple and does not require a CNS
solver, it seems to be not accurate enough.
\end{remark}

   \subsection{A 3D example}
   \label{subsec:3D}

\begin{figure}[h]
    \centering
    \includegraphics[width=0.5\textwidth]{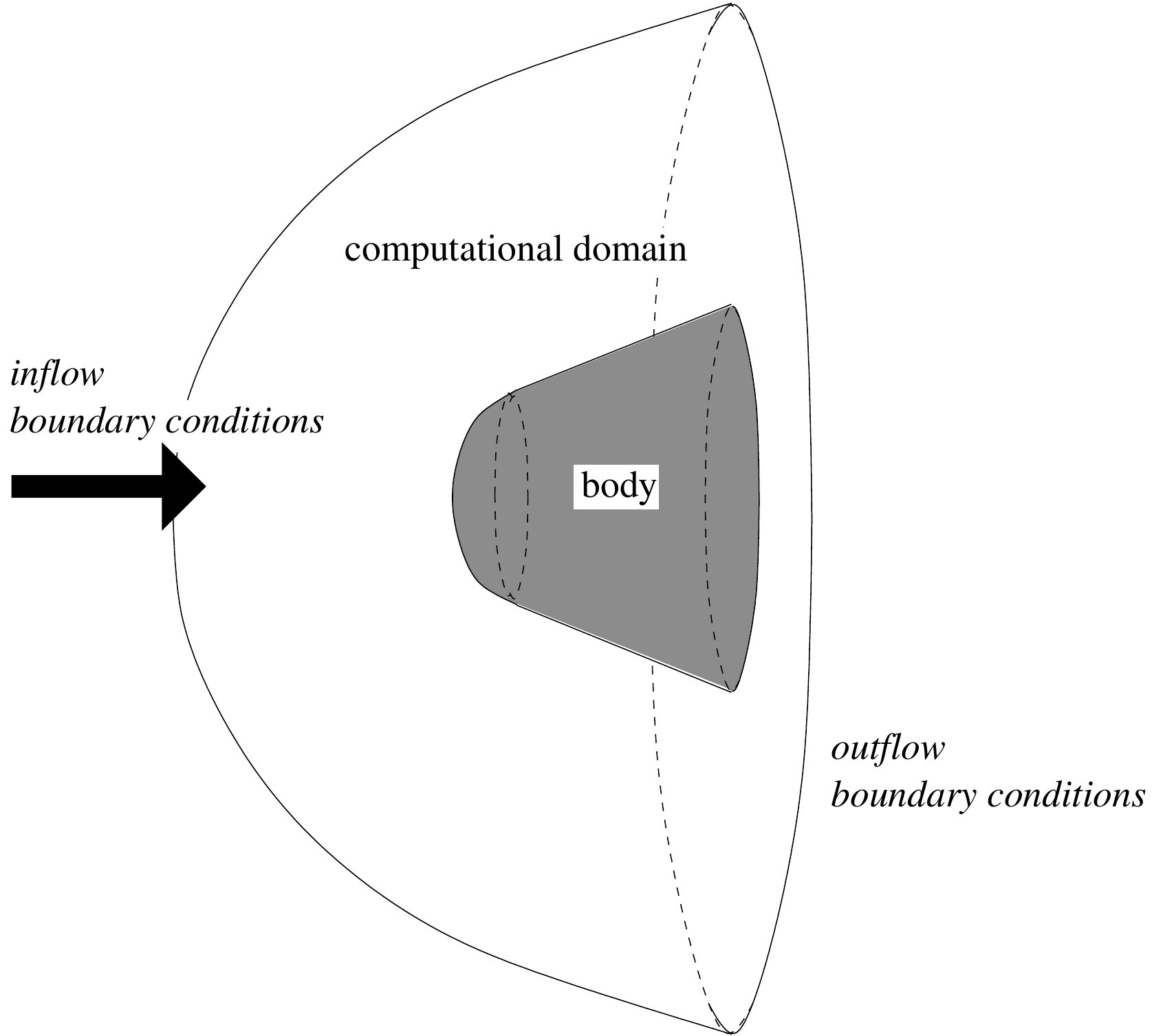}
    \caption{\modif{3D flow around a cone: geometry and
      computational domain. The downstream flow is not simulated.}}
    \label{fig:geom_3D}
  \end{figure}

   Here, we consider the flow around a 3D object composed of cone with
   a spherical nose, with no incidence \modif{(see
     figure~\ref{fig:geom_3D})}. The physical parameters are the same
   as in section~\ref{subsec:res_axi}. The space mesh uses $50\, 000$
   cells, with 50 in the direction normal to the
   surface. \modif{Again, the downstream flow is ignored, and outflow
     boundary conditions are used at the boundary of the right side of
     the domain. }

A CNS pre-simulation and the use of
formula~(\ref{eq-bounds_steps_axi}) with parameter $c=4$ and $a=2$
 give us a uniform Cartesian grid with $65\times33\times 33=70\, 785$
 discrete velocities. This is
 much too large to make the simulation (the memory storage itself is
 huge).  Then we only do the simulation with the AMR grid (${\mathbb
   P}_0$ version) generated by our algorithm which has only $2\, 956$
 points.

The gain factor in memory storage is 22.  The steady state is
 reached in $3 \, 424$ iterations and $3\, 566 s$ with 400 domains
 and 4 threads per domain. 

 To estimate the gain in CPU that would be obtained if the uniform
 Cartesian grid could be used, we use the following remark: on the
 previous tests (2D plane and axisymmetric), the number of iterations
 to reach the steady state is almost the same for the uniform and the AMR
 velocity grids. Then we make a single iteration with the uniform
 Cartesian grid (same domain decomposition) and we multiply the CPU
 time obtained for this iteration by the number of iterations
 required with the AMR grid.  We find a CPU time which is 24 times as
 large as with the AMR grid, which is quite important. This CPU gain
 is larger that the gain observed for the memory storage.

\begin{figure}[h]
    \centering
   \includegraphics[height=0.4\textheight]{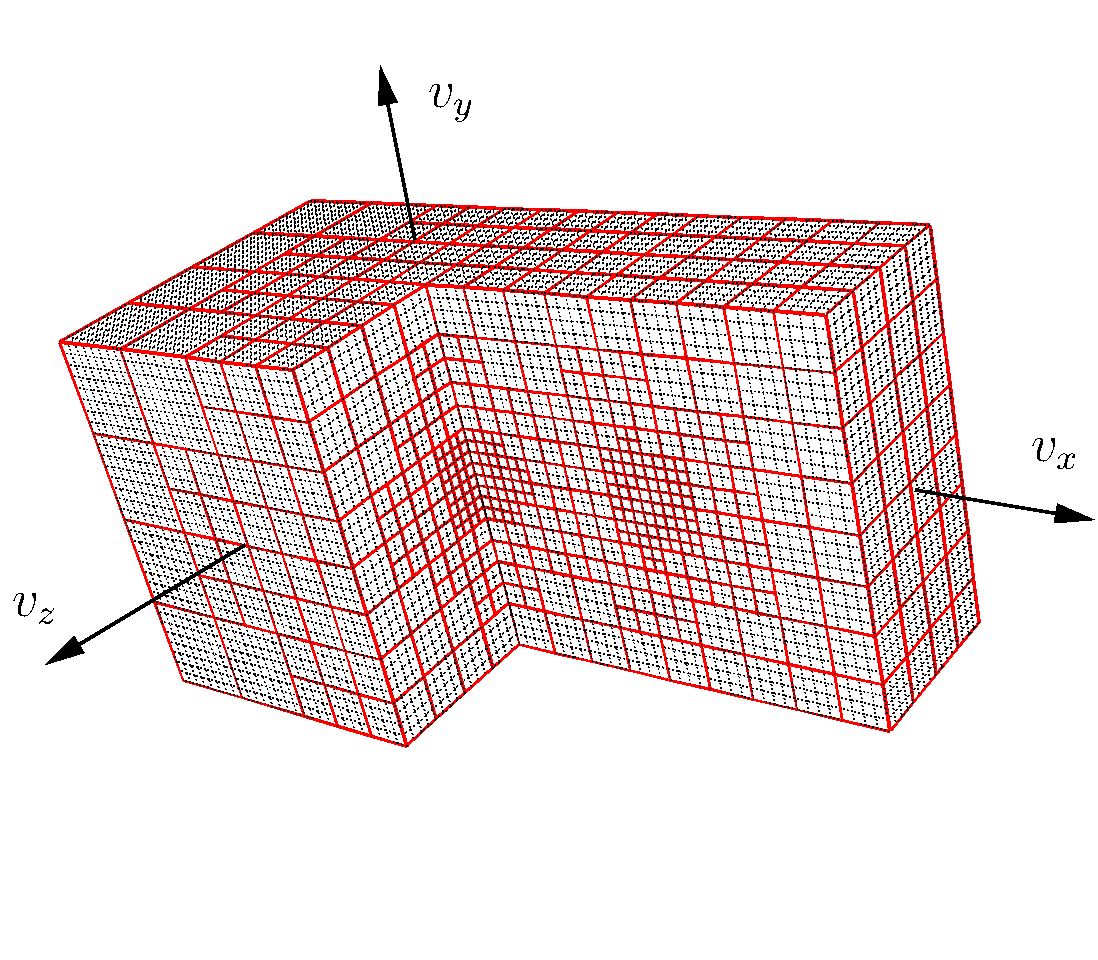}
    \includegraphics[height=0.4\textheight]{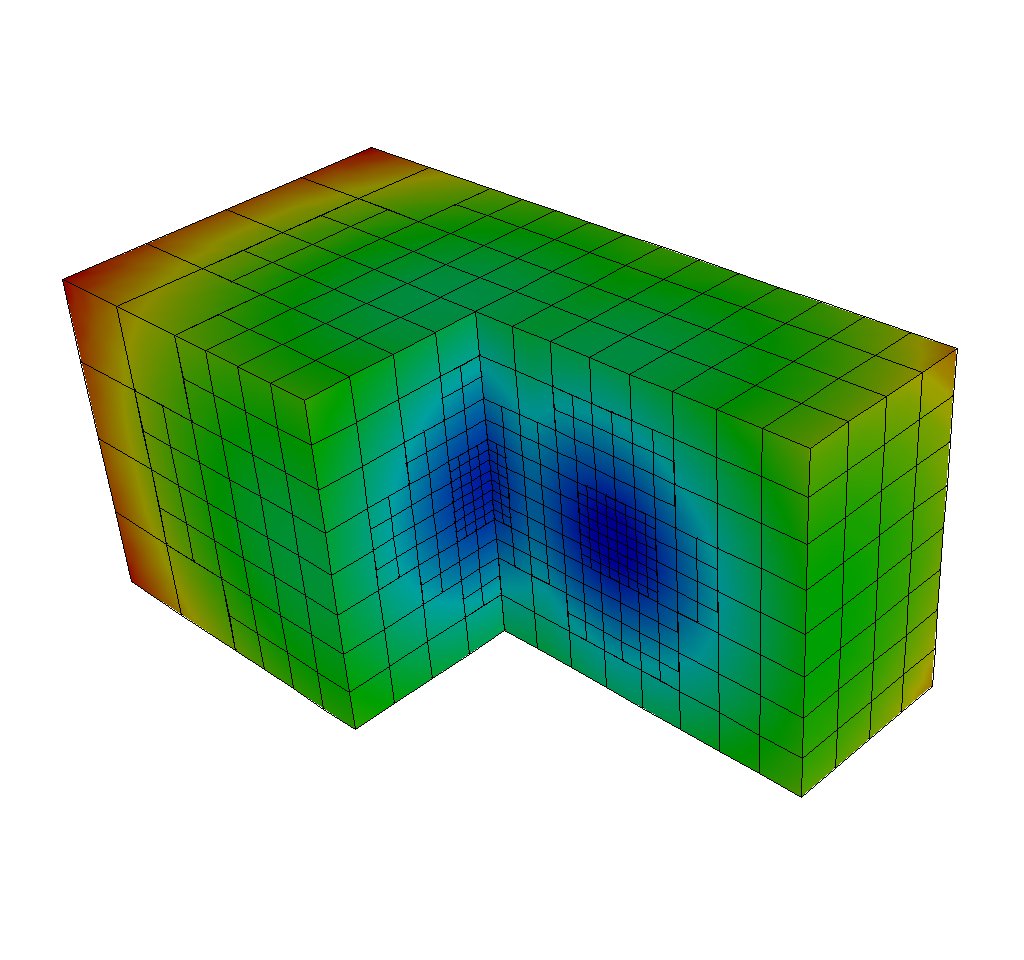}
    \caption{
      \modif{3D flow around a cone. Top: the 3D Cartesian velocity
        grid (dashed line) and AMR velocity grid (solid) with the
        Cartesian coordinate system $(v_x,v_y,v_z)$. In order to see
        the refined zones around $(0,0,0)$ and $(u_{upstream},0,0)$
        that are inside the grid, a part of the grid has been
        hidden. Bottom: The same AMR velocity grid with the contours
        of the corresponding support function based on the CNS
        pre-simulation.}  }
    \label{fig:3D_grids}
  \end{figure}

\begin{figure}[h]
    \centering
    \includegraphics[width=0.6\textwidth]{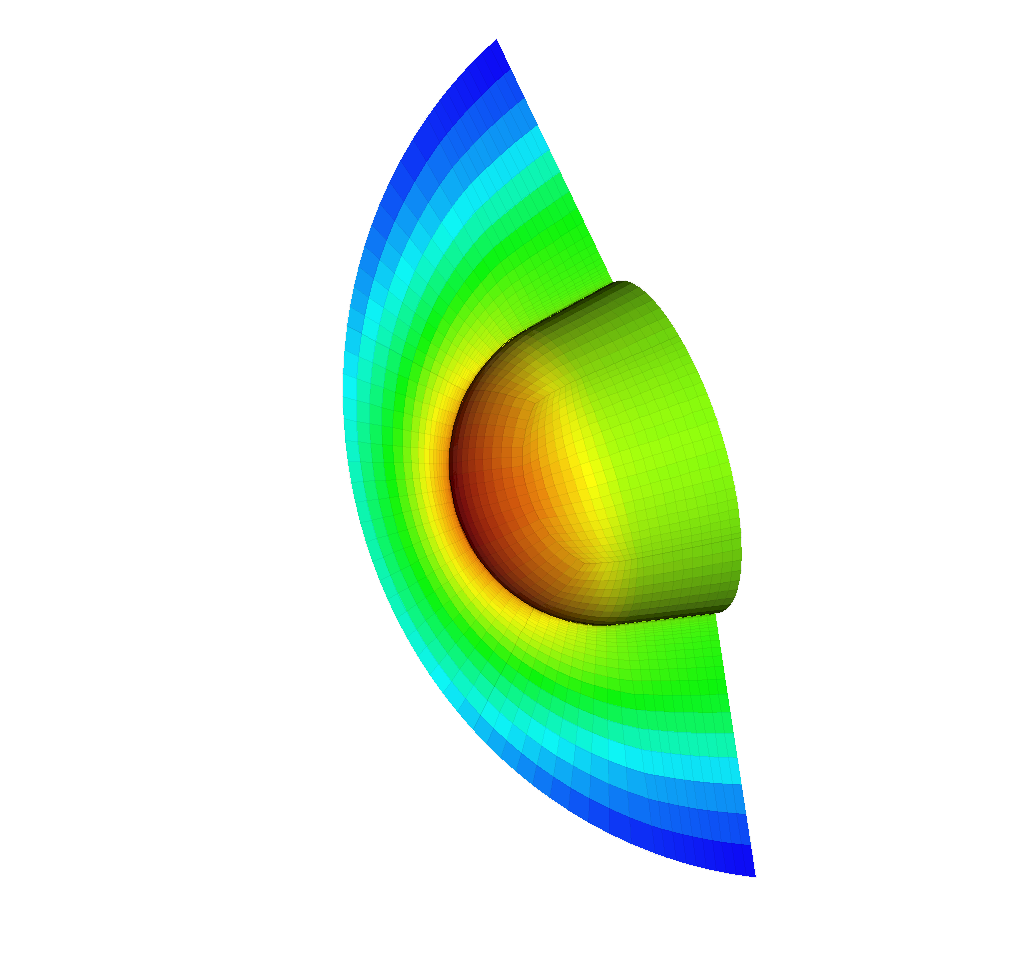}
    \caption{3D flow around a cone: pressure field computed with the
      kinetic code and the AMR velocity grid. \modif{The pressure is shown on
    the surface of the solid, as well as in a vertical plane $(x,y,0)$ that crosses the cone
along its middle part}.}
    \label{fig:3D_pressure}
  \end{figure}

\begin{figure}[h]
    \centering
     \includegraphics[width=\textwidth]{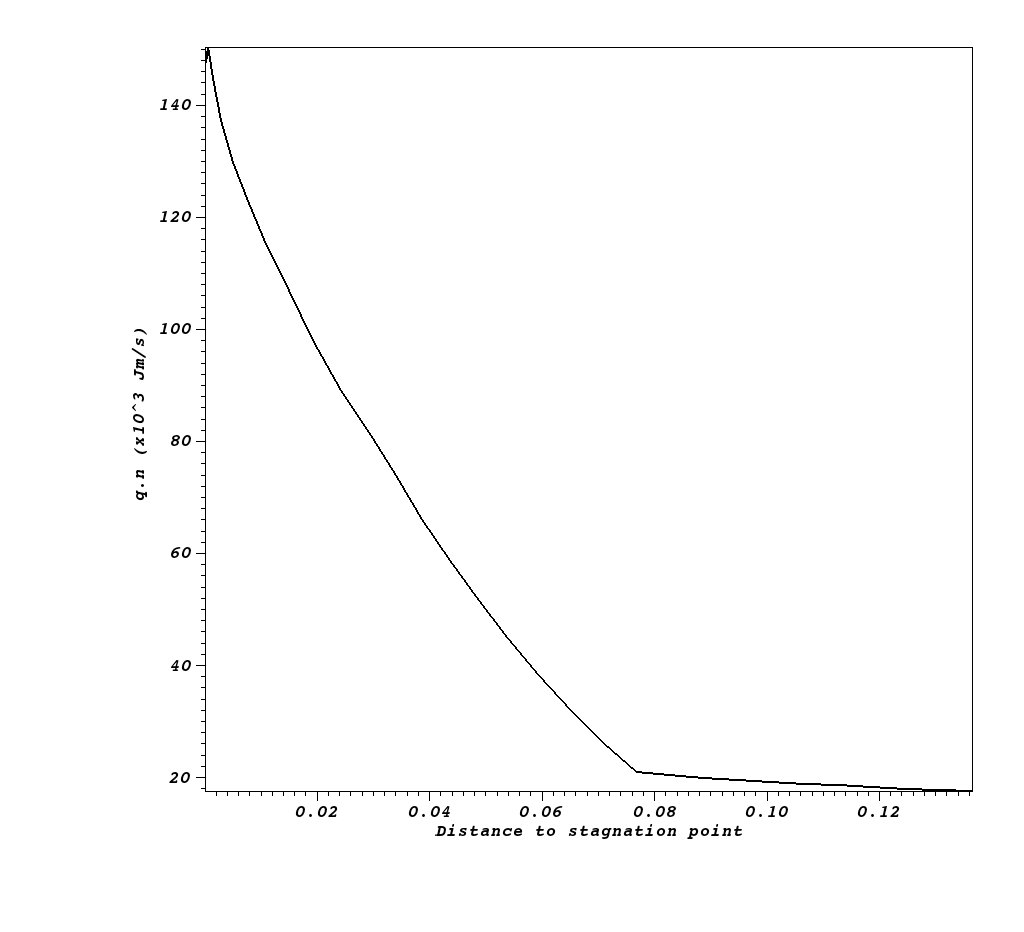}
    \caption{3D flow around a cone: normal heat flux along the
      boundary (computed with the
      kinetic code and the AMR velocity grid).}
    \label{fig:3D_flux}
  \end{figure}

In figure~\ref{fig:3D_grids}, we show the uniform Cartesian grid and
the AMR grid. In figure~\ref{fig:3D_pressure}, we show the pressure
field computed by the code with the AMR grid, and the heat flux along
the surface is shown in figure~\ref{fig:3D_flux}.

   \subsection{Initialization with the CNS results}
   \label{subsec:init}

   Since we use a CNS pre-simulation to build our discrete velocity
   grid, it is interesting to use it also to initialize the BGK
   computation, what was done in the previous simulations. It helps to
   reach the steady state more rapidly, as it is shown in the
   following example.

We take the same test case as in
   section~\ref{subsec:res_axi}, with the uniform velocity grid. In a
   first computation, the distribution is initialized with the
   (Maxwellian) upstream flow (this is the standard initialization),
   and the scheme converges to the steady state in $3\, 549$
   iterations. Then we do the same computation, but now the initial
   state is the Maxwellian distribution computed with the macroscopic
   variables that are given by the steady CNS solution. Then the
   steady state is reached in $1\, 822$ iterations only.

A similar gain (half as many iterations as with the standard
initialization) is obtained in all our test cases. 

\section{Conclusion and perspectives}
\label{sec:conclusion}
In this paper, we have proposed a method to refine and coarsen a
global velocity grid of a discrete velocity approximation of a kinetic
equation. It is based on a criterion called ``the support function''
that links the local size of the velocity grid to the macroscopic
temperature an velocity of the flow. Our algorithm uses the
macroscopic fields given by a compressible Navier-Stokes
pre-simulation to automatically generates an optimal velocity grid.

This approach has been tested in a 3D computational code (and 2D plane
and axisymmetric versions) which uses the BGK model of rarefied gas
dynamics. Typical test cases in hypersonic re-entry problems (for
simplified geometries) have been used.  We observed that using our
algorithm allows important gains both in CPU time and memory storage,
up to a factor 24.  This allows to make
simulations that are hardly possible with standard grids, even on
super computers. 

\modif{Note that this work might be extended to unsteady flow
  simulations: we could indeed modify the velocity grid at different
  time steps according to our refinement algorithm. This might also be
  used during the iterations of our iterative solver for a steady
  simulation. However, this would require to interpolate the solution
  between two successive grids, and this might lead to an important
  increase of the CPU time. This clearly requires more investigations
  to get an efficient method.}

An obvious, and straightforward, extension of this work will be the use
of modified BGK models (like ES or Shakhov models), in order to have a
model with a correct Prandtl number. 

Moreover, we will try to investigate non isotropic AMR velocity grids
by taking into account translational kinetic temperatures in different
directions, that can be obtained by the stress tensor computed by the
CNS simulation. \modif{It could also be interesting to refine the grid
  so as to capture the discontinuity of the distribution function in
  the velocity space, but this will require a rather different approach.}

Our main short term project is to improve the physical accuracy of our
code by implementing a simple BGK-like model for multi-species
reactive flows. An application of our method to the full Boltzmann
collision operator is less obvious but might also be explored.

\paragraph{Acknowledgment.} Experiments presented in this paper were
carried out using the PLAFRIM experimental testbed, being developed
under the Inria PlaFRIM development action with support from LABRI and
IMB and other entities: Conseil R\'egional d'Aquitaine, FeDER,
Universit\'e de Bordeaux and CNRS (see
https://plafrim.bordeaux.inria.fr/).

 \bibliographystyle{plain}

 \bibliography{biblio}

\appendix
\section{Overview of the 3D kinetic code}
\label{sec:code}

   \subsection{The linearized implicit scheme}
Our code is an extension of the code presented
in~\cite{luc_jcp} to 3D polyatomic flows. It is based on the following
reduced BGK model
 \begin{align*}
  & \dt f + \vit \cdot \nabla_\x f = \frac{1}{\tau}(M(U)-f) \\
  & \dt g + \vit \cdot \nabla_\x g = \frac{1}{\tau}(\frac{\delta}{2}RTM(U)-g),
 \end{align*}
where $U=\cint{mf + e^{(5)} g}=(\rho,\rho\um,E=\demi\rho|\um|^2+\frac{3+\delta}{2}\rho R T)$
is the vector of macroscopic mass, momentum, and energy
density. Here, we use the standard notation
$\cint{.}=\int_{\R^3}.\,d\vit$ for any vector valued function of
$\vit$, $m(v)=(1,\vit,\demi|\vit|^2)$ is the vector of collisional
invariants, and $e^{(5)}=(0,0,0,0,1)$. Moreover, $M(U)$ is the standard
Maxwellian distribution defined through the density, velocity, and
temperature corresponding to the vector $U$ above
(see~(\ref{eq-maxw})). 

This model comes from the reduction of a BGK model for the full
distribution function $F(t,\x,\vit,I)$, where $I$ is the internal
energy variable, and $\delta$ is the number of internal degrees of
freedom (see~\cite{HH_1968} for the first use of this technique
and~\cite{ALPP} for an application to polyatomic gases). Consequently,
it accounts for any number of internal degrees of freedom. For
instance, a diatomic gas can be described with $\delta=2$.

This model is first discretized with respect to the velocity
variable. We follow the approximation of~\cite{luc_jcp} and its
extension to polyatomic gases~\cite{luc_dubroca}. We assume we have a
velocity grid $\lbrace \vit_q,q=0:q_{max}  \rbrace$ (like the grids
described in this paper). The continuous distributions $f$ and $g$ are
then replaced by their approximations at each point $\vit_q$, and we
get the following discrete velocity BGK system
 \begin{align*}
  & \dt f_q + \vit_q \cdot \nabla_\x f_q = \frac{1}{\tau}(M_q(U)-f_q) \\
  & \dt g_q + \vit_q \cdot \nabla_\x g_q = \frac{1}{\tau}(N_q(U)-g_q),
 \end{align*}
where $(M_q(U),N_q(U))$ is an approximation of
$(M(U),\frac{\delta}{2}RTM(U))$ that has to be defined. As explained
in this paper, we assume we have quadrature weights
$\lbrace\omega_q\rbrace$ corresponding to our discrete velocity grid,
so that the moment vector of the discrete distributions is
\begin{equation*}
U=\sum_{q=0}^{q_{max}}(m(\vit_q)f_q+e^{(5)}g_q)\omega_q. 
\end{equation*}
As proposed
in~\cite{luc_dubroca}, the discrete equilibrium $(M_q(U),N_q(U))$ is
constructed so that it has the same moments as $(f,g)$, that is to say
\begin{equation}\label{eq-discM} 
  \sum_{q=0}^{q_{max}}(m(\vit_q)M_q(U)+e^{(5)}N_q(U))\omega_q = U.
\end{equation}
In our code, $(M_q(U),N_q(U))$ is determined through the entropic
variable $\alpha$ such that $M_q(U)=\exp(\alpha\cdot m(\vit_q))$ and
$N_q(U)=\frac{\delta}{2}\frac{1}{-\alpha_5}\exp(\alpha\cdot
m(\vit_q))$, by solving~(\ref{eq-discM}) by a Newton algorithm. We
mention that the computational cost of this algorithm can be
significantly reduced by using a nice idea due to
Titarev~\cite{T_cicp}. This optimization will be used in a future
version of our code.

The discrete velocity BGK system is then discretized by a finite
volume scheme on a multiblock curvilinear 3D mesh of hexahedral cells
$\Omega_{ijk}$, with indices $i,j,k=1$ to $i_{max},j_{max},k_{max}$,
respectively. Denoting by $\fnijkq$ an approximation of the average of
$f$ at time $t_n$ on a cell $\Omega_{ijk}$ at the discrete velocity
$\vit_q$, our scheme reads, in its implicit version,
\begin{align*}
&   \frac{\fnpunijkq-\fnijkq}{\Dt} + \left(\vit_q\cdot
    \nabla_{\x}f^{n+1}_q\right)_{ijk} =
  \frac{1}{\tau^{n+1}_{ijk}}(M_q(\Unpunijk)-\fnpunijkq) , \\
&   \frac{\gnpunijkq-\gnijkq}{\Dt} + \left(\vit_q\cdot
    \nabla_{\x}g^{n+1}_q\right)_{ijk} = \frac{1}{\tau^{n+1}_{ijk}}(N_q(\Unpunijk)-\gnpunijkq) ,
\end{align*}
where 
\begin{equation*}
\Unpunijk=\sum_{q=0}^{q_{max}}(m(\vit_q)\fnpunijkq+e^{(5)}\gnpunijkq)\omega_q. 
\end{equation*}
The discrete divergence $\left(\vit_q\cdot
  \nabla_{\x}f_q\right)^{n+1}_{ijk}$ is given by the following second
order upwind approximation (with the Yee limiter~\cite{yee}):
\begin{equation}\label{eq-divergence}
  \left(\vit_q\cdot \nabla_{\x}f^{n+1}_q\right)_{ijk} =
\frac{1}{|\Omega_{i,j,k}|}\left( 
\Bigl(\Phi_{i+\demi,j,k} - \Phi_{i-\demi,j,k}\Bigr)
+ \Bigl(\Phi_{i,j+\demi,k} - \Phi_{i,j-\demi,k}\Bigr)
+ \Bigl(\Phi_{i,j,k+\demi} - \Phi_{i,j,k-\demi}\Bigr)
   \right),
\end{equation}
where 
\begin{equation*}
\begin{split}
  \Phi_{i+\demi,j,k}& =(\vit_q\cdot\nu_{i+\demi,j,k})^+\fnpunijkq +
  (\vit_q\cdot\nu_{i+\demi,j,k})^- f^{n+1}_{i+1,j,k,q} \\
 & \quad +
  |\vit_q\cdot\nu_{i+\demi,j,k}|\, {\rm minmod}(\Delta^{n+1}_{i-\demi},\Delta^{n+1}_{i+\demi},\Delta^{n+1}_{i+\frac{3}{2}})
\end{split}
\end{equation*}
is the second order numerical flux across the face between
$\Omega_{i,j,k}$ and $\Omega_{i+1,j,k}$, and $\nu_{i+\demi,j,k}$ is
the normal vector to this face directed from $\Omega_{i,j,k}$ to
$\Omega_{i+1,j,k}$ while its norm is equal to the area of the face. In
the minmod limiter function, we use the notation
$\Delta^{n+1}_{i+\demi}=f^{n+1}_
{i,+1,j,k,q}-\fnpunijkq$. Finally, we use the standard notation
$a^\pm=(a\pm|a|)/2$ for every number $a$. The numerical
fluxes across the other faces are defined accordingly. For the sequel,
it is useful to denote by $\left(\vit_q\cdot
  \nabla_{\x}f^{n+1}_q\right)^{1st\ order}_{ijk}$ the corresponding linear
first order upwind discretization (that is to say, defined
by~(\ref{eq-divergence}) where the minmod term is set to 0).

The advantage of the time implicit approximation is that it ensures
unconditional stability, which allows us to take large time steps, and
hence to get rapid convergence to the steady state. Of course, this scheme is too expensive, since it requires to solve a
non linear system at each time iteration. Therefore, we rather use the
following linearization. First, the Maxwellian is linearized by a
first order Taylor expansion:
\begin{equation*}
  M_q(\Unpunijk)=M_q(\Unijk+(\Unpunijk-\Unijk)) \approx M_q(\Unijk)
  + \partial_U M_q(\Unijk)(\Unpunijk-\Unijk),
\end{equation*}
and the same for $N_q(\Unpunijk)$ (see
section~\ref{appendix:jacobian} for explicit expressions of the
Jacobian matrices). Then the discrete divergence, which is not
differentiable due to the limiter, is linearized by using the
corresponding first order upwind approximation. This is sometimes
called a ``frozen coefficient technique'', which gives 
\begin{equation*}
\begin{split}
  \left(\vit_q\cdot \nabla_{\x}f_q\right)^{n+1}_{ijk} & =
  \left(\vit_q\cdot \nabla_{\x}(f^n_q+(f^{n+1}_q-f^n_q))\right)_{ijk}
  \\
  & \approx \left(\vit_q\cdot \nabla_{\x}f^n_q\right)_{ijk}
  +\left(\vit_q\cdot \nabla_{\x}(f^{n+1}_q-f^n_q)\right)^{1st\
    order}_{ijk}
\end{split}
\end{equation*}
Then, denoting by $\delta f_{i,j,k,q}=\fnpunijkq-\fnijkq$ (same
notation for $g$), and by $\delta U_{i,j,k}$ the moments of $(\delta
f_{i,j,k},\delta g_{i,j,k})$, our scheme reads in the following form:
\begin{align}
&  \frac{\delta f_{i,j,k,q}}{\Dt} + \left(\vit_q\cdot
    \nabla_{\x}\delta f_q\right)^{1st\ order}_{ijk} -
  \frac{1}{\tau^n_{i,j,k}}\partial_U M_q(\Unijk)(\delta U_{i,j,k}) =
  RHSf^n_{i,j,k}\label{eq-dformf}   \\
&  \frac{\delta g_{i,j,k,q}}{\Dt} + \left(\vit_q\cdot
    \nabla_{\x}\delta g_q\right)^{1st\ order}_{ijk} -
  \frac{1}{\tau^n_{i,j,k}}\partial_U N_q(\Unijk)(\delta U_{i,j,k}) =
  RHSg^n_{i,j,k},\label{eq-dformg} 
\end{align}
where the right-hand sides are given by
\begin{align*}
  RHSf^n_{i,j,k} & = 
- \left(\vit_q\cdot \nabla_{\x}f^{n}_q\right)_{ijk}
+  \frac{1}{\tau^{n}_{ijk}}(M_q(\Unijk)-\fnijkq) \\
  RHSg^n_{i,j,k} & = 
- \left(\vit_q\cdot \nabla_{\x}g^{n}_q\right)_{ijk}
+  \frac{1}{\tau^{n}_{ijk}}(N_q(\Unijk)-\gnijkq).
\end{align*}

If our scheme converges to steady state, then the right-hand side is
zero, and we get a second order discrete steady solution.

\modif{
   \subsection{Numerical boundary conditions}
   \label{subsec:num_bc}

   Numerically, the boundary conditions are implemented by the
   standard ghost cell technique, which is used as follows. When an
   index $(i,j,k)$ corresponds to a cell located at the boundaries of
   the domain, there appear unknown values in the numerical fluxes
   like $f^n_{i,0,k,q}$ and $f^n_{i,j_{max}+1,k,q}$, for the cells
   $(i,1,k)$ and $(i,j_{max},k)$ for instance (see~(\ref{eq-divergence})). Corresponding cells $(i,0,k)$,
   $(i,j_{max}+1,k)$, etc. are called ghost-cells.  These values are
   classically defined according to the boundary conditions that are specified
   for the problem. Here we use several kinds of boundary conditions:
   solid wall interactions, inflow and outflow boundary conditions at
   artificial boundaries, as well as symmetry boundary conditions
   along symmetry planes and symmetry axes.

For the diffuse reflection, the incident molecules in a boundary cell of
index $(i,1,k)$ are supposed to be re-emitted by the wall from a
ghost cell of index $(i,0,k)$. This cell is the mirror cell of
$(i,1,k)$ with respect to the wall. The diffuse
reflection~(\ref{eq-BCdiff})--(\ref{eq-sigmadiff}) is then
modeled by
\begin{equation} \label{paper-eq-BCdiff_num} 
(f^n_{i,0,k,q},g^n_{i,0,k,q}) =
  \sigma_{i,1,k} \, (M^{wall}_q,N_q^{wall}), \quad {\bf v}_q
  \cdot {\rm n}_{i,1,k} >0,
\end{equation}
where $\sigma_{i,1,k}$ is determined so as to avoid a mass flux across
the wall, that is to say between cells $(i,0,k)$ and $(i,1,k)$. In
this relation, $(M^{wall}_q,N_q^{wall})$ is a discrete conservative
approximation of the wall Maxwellians $(1,\frac{\delta}{2}
RT_{wall})\frac{1}{(2\pi RT_{wall})^{3/2}}\exp(-|{\bf
  v}|^2/2RT_{wall})$.  Relation~(\ref{eq-sigmadiff}) gives
$$
\sigma_{i,1,k} = - \frac{\sum_{{\bf v}_q \cdot n_{i,1,k}<0} {\bf v}_q \cdot
  {\rm n}_{i,1,k} \, f^n_{i,1,k,q} \omega_q} {\sum_{{\bf v}_q \cdot
    n_{i,1,k}>0} {\bf v}_q \cdot n_{i,1,k} \, M^{wall}_q \omega_q}.
$$

For the inflow boundary condition, for instance at a boundary cell
$(i,j_{max}+1,k)$, we simply set the ghost cell value to the upstream
Maxwellian distributions $(M^{upstream},N^{upstream})$ (defined
through the upstream values of $\rho$, $u$, and $T$):
$$
(f^n_{i,j_{max}+1,k,q},g^n_{i,j_{max}+1,k,q}) =(M_{q}^{upstream},N_{q}^{upstream}).
$$

For the outflow boundary condition, we set the ghost cell value to
the value of the corresponding boundary cell:
$$
(f^n_{i,j_{max}+1,k,q},g^n_{i,j_{max}+1,k,q}) =
(f^n_{i,j_{max},k,q},g^n_{i,j_{max},k,q}).
$$

Finally, for a cell $(1,j,k)$ in a symmetry plane (for instance the
plane $(0,y,z)$) we use the symmetry of the distribution functions to
set
$$
(f^n_{0,j,k,q},g^n_{0,j,k,q}) =
(f^n_{1,j,k,q'},g^n_{1,j,k,q'}),
$$
where $q'$ is such that ${\bf v}_{q'}$ is the symmetric of ${\bf v}_q$
with respect to the symmetry plane. In this case, our AMR velocity
grid is constructed such that it is also symmetric with respect to
this plane, which gives ${\bf
  v}_{q'}=(-v_{x,q},v_{y,q},v_{z,q})$. This grid is obtained in two
steps: first, a part on one side of the symmetry plane is
obtained by using our algorithms (support function and AMR grid
generation), and then this part is symmetrized to obtain the part on
the other side.

For the second order numerical flux, we also need the values of a
second layer of ghost cells with indices like $(-1,j,k)$ or
$(i_{max}+2,j,k)$, etc. For these values, we simply copy the value of the
corresponding ghost cell in the first layer, that is to say
$f^n_{-1,j,k,q}=f^n_{0,j,k,q}$ for instance. This treatment makes the
accuracy of our scheme reduce to first order at the boundary (since
this makes the flux limiters equal to 0). It may
be more relevant to use extrapolation techniques, which is the
subject of a work in progress. Only the boundary
condition on a symmetry plane is treated differently: here, we use
the symmetry of the distributions to set
$$
(f^n_{-1,j,k,q},g^n_{-1,j,k,q}) =
(f^n_{2,j,k,q'},g^n_{2,j,k,q'}),
$$
where $v_{q'}$ has been defined above.
}

   \subsection{Linear solver}

At each time iteration, our linearized implicit scheme requires to
solve the large linear system~(\ref{eq-dformf}--\ref{eq-dformg}). It is
therefore interesting to write it in the following matrix form:
\begin{equation*}
  \Bigl( \frac{I}{\Dt} + T + R^n  \Bigr) \delta F = RHS^n,
\end{equation*}
where $\delta F=(\delta f_{i,j,k,q} , \delta g_{i,j,k,q})$ is a large
vector that stores all the unknowns, $I$ is the unit matrix, $T$ is
a matrix such that 
\begin{equation}\label{eq-matrixT}
  (T\delta F)_{i,j,k,q}=\left( 
\Bigl(\vit_q\cdot \nabla_{\x}\delta f_q\Bigr)^{1st\ order}_{ijk}   
, \Bigl(\vit_q\cdot \nabla_{\x}\delta g_q\Bigr)^{1st\ order}_{ijk}   
\right),
\end{equation}
$R^n$ is the relaxation matrix such that
\begin{equation}\label{eq-matrixRn} 
  (R^n\delta F)_{i,j,k,q}= \left( - \frac{1}{\tau^n_{i,j,k}}\partial_U  M_q(\Unijk)(\delta U_{i,j,k})
,   - \frac{1}{\tau^n_{i,j,k}}\partial_U  N_q(\Unijk)(\delta U_{i,j,k})
\right),
\end{equation}
and $RHS^n= \left( RHSf^n_{i,j,k}, RHSg^n_{i,j,k}\right)$.

For simplicity, we use explicit boundary conditions, which means $\delta
F_{i,j,k,q}=0$ in ghost cells. This implies that $T$ has a simple
block structure, which is used in the following. 

The algorithm used in our code is based on a coupling between the
iterative Jacobi and Gauss-Seidel methods. First, the relaxation
matrix $R^n$ is splitted into its diagonal $\Delta^n$ and its
off-diagonal $-E^n$, so that we get the following system (this the
Jacobi iteration): 
\begin{equation}\label{eq-jacobi} 
  \Bigl( \frac{I}{\Dt} + T + \Delta^n  \Bigr) \delta F = RHS^n + E^n
  \delta F.
\end{equation}
Note that $E^n$ is very sparse: the product $E^n\delta F$ is local in
space and can be written $[E^n_{i,j,k}\delta F_{i,j,k}]_q$ (see
section~\ref{appendix:jacobian}). The left-hand side of this system
has a three level tridiagonal block structure. We solve it by using a
line Gauss-Seidel iteration which is described below. 

With curvilinear grids, in many cases, and in particular for re-entry
problems, the flow is aligned with a mesh (generally aligned
with solid boundaries), which means that the largest variations occur
along the orthogonal direction, say the direction of index $i$, for
instance. The idea is to use an ``exact'' inversion of $T$ along this
direction. This is done with a sweeping strategy sometimes called
``Gauss-Seidel line iteration''. First, we rewrite
system~(\ref{eq-jacobi}) pointwise, as follows:
\begin{equation*}
\begin{split}
& \left(\frac{1}{\Dt}+A+B+C+\Delta^n_{i,j,k,q}\right)\delta F_{i,j,k,q}
 + A^-\delta F_{i+1,j,k,q}  + A^+\delta F_{i-1,j,k,q} \\
& + B^-\delta F_{i,j+1,k,q}  + B^+\delta F_{i,j-1,k,q}
+ C^-\delta F_{i,j,k+1,q}  + C^+\delta F_{i,j,k-1,q} 
=RHS^n_{i,j,k,q} +[E^n_{i,j,k}\delta F_{i,j,k}]_q,
\end{split}
\end{equation*}
where coefficients $A,B,C,A^\pm,B^\pm,C^\pm$ are standard notations
for line-Gauss-Seidel method and can be easily derived
from~(\ref{eq-jacobi}), (\ref{eq-matrixT}),
and~(\ref{eq-divergence}). Then the linear solver reads as shown in
Algorithm~\ref{alg:JGS}. \\
\begin{algorithm}       
\caption{Jacobi/line-Gauss-Seidel algorithm}
\label{alg:JGS}
\begin{algorithmic}
\State $\delta F^{(0)}=0$ \Comment{initialization}
  \For{$p$ from 0 to $P$}   \Comment{iterations of the solver}
          \For{$q$ from 0 to $q_{max}$}  \Comment{loop over the
            discrete velocities (Jacobi loop)}
\State

           \Comment{one sweep in $j$ direction (forward substitution)}
            \For{$k$ from 1 to $k_{max}$}  
                 \For{$j$ from 1 to $j_{max}$} 
\State 
             solve exactly the tridiagonal system 
\begin{equation*}
\begin{split}
& \left(\frac{1}{\Dt}+A+B+\Delta^n_{i,j,k,q}\right)\delta F^{(p+\demi)}_{i,j,k,q}
 + A^-\delta F^{(p+\demi)}_{i+1,j,k,q}  + A^+\delta F^{(p+\demi)}_{i-1,j,k,q} \\
& \qquad =RHS^n_{i,j,k,q} +[E^n_{i,j,k}\delta F^{(p)}_{i,j,k}]_q
-  B^+\delta F^{(p+\demi)}_{i,j-1,k,q} - B^-\delta F^{(p)}_{i,j+1,k,q}
\end{split}
\end{equation*}
                  \EndFor

                \EndFor
 \State 
 
              \Comment{one sweep in $k$ direction (forward  substitution)}  
            \For{$j$ from 1 to $j_{max}$}  
                 \For{$k$ from 1 to $k_{max}$} 
\State 
             solve exactly the tridiagonal system 
\begin{equation*}
\begin{split}
& \left(\frac{1}{\Dt}+A+C+\Delta^n_{i,j,k,q}\right)\delta F^{(p+1)}_{i,j,k,q}
 + A^-\delta F^{(p+1)}_{i+1,j,k,q}  + A^+\delta F^{(p+1)}_{i-1,j,k,q} \\
& \qquad =RHS^n_{i,j,k,q} +[E^n_{i,j,k}\delta F^{(p)}_{i,j,k}]_q
-  C^+\delta F^{(p+1)}_{i,j,k-1,q} - C^-\delta F^{(p+\demi)}_{i,j,k+1,q}
\end{split}
\end{equation*}
                  \EndFor
           \EndFor
      \EndFor
  \State 
       \State compute the moments of $\delta F^{(p+1)}$ \Comment{for the
           action of $\Delta^n_{i,j,k,q}$ and $E^n_{i,j,k}$ onto
           $\delta F^{(p+1)}$}
       \State
    \EndFor
    \State     set $\delta F = \delta F^{(P+1)}$
  \end{algorithmic}
\end{algorithm}

\begin{remark}
  \begin{enumerate}
  \item We observe that we get a fast convergence to steady state by
    using a few step of this linear solver (say $P=2$ or 3). This is
    due to the ``exact'' inversion along the direction of largest
    variation of the flow.
    \item In practice, we add one back substitution step right after
      each forward substitution step. 
    \item This solver is a straightforward extension of the linear
      solver used in~\cite{luc_jcp}. More sophisticated versions of
      this solver could be tested.
\end{enumerate}
\end{remark}

   \subsection{Jacobian matrices and the relaxation matrix $R^n$}
   \label{appendix:jacobian}

Elementary calculus gives the following formula:
\begin{align*}
& \partial_{U} M_q(U)    = M_q(U) m(\vit_q) A(U)^{-1} \\
& \partial_{U} N_q(U)    = N_q(U) (m(\vit_q)-\frac{1}{\alpha_5(U)}e^{(5)}\,) A(U)^{-1},
\end{align*}
where $A(U)$ is the following $5\times 5$ matrix
\begin{equation*}
  A(U)=\sum_{q=0}^{q_{max}} \left(m(\vit_q)^Tm(\vit_q)M_q(U) +
    {e^{(5)}}^T\Bigl(m(\vit_q)-\frac{1}{\alpha_5(U)}e^{(5)}\Bigr)N_q(U)\right)\, \omega_q.
\end{equation*}
Consequently, by using the definition~(\ref{eq-matrixRn}) of the
relaxation matrix $R^n$, a careful algebra shows that its diagonal elements
$\Delta^n_{i,j,k,q}$ can be separated into the following two blocks $(\Delta^{n,f}_{i,j,k,q},\Delta^{n,g}_{i,j,k,q})$:
\begin{align*}
& \Delta^{n,f}_{i,j,k,q}=\frac{1}{\tau^n_{i,j,k}}\left(M_q[\Unijk] m(\vit_q)
A(\Unijk)^{-1} m(\vit_q)^T\omega_q -1\right),\\
& \Delta^{n,g}_{i,j,k,q}=\frac{1}{\tau^n_{i,j,k}}\left(N_q[\Unijk] \Bigl(m(\vit_q)-\frac{1}{\alpha_5(\Unijk)}e^{(5)}\Bigr) A(\Unijk)^{-1} {e^{(5)}}^T\omega_q-1\right).
\end{align*}
Therefore, the product of the off-diagonal part $-E^n$ of $R^n$ with
$\delta F$ (as used in algorithm~\ref{alg:JGS}) is simply
\begin{equation*}
\begin{split}
&  -[E^n_{i,j,k}\delta F_{i,j,k}]_q=\left( 
 - \frac{1}{\tau^n_{i,j,k}}\partial_U  M_q(\Unijk)(\delta U_{i,j,k})
- \Delta^{n,f}_{i,j,k,q}\delta f_{i,j,k,q}
 \right.  \\
& \qquad\qquad\qquad\qquad\qquad, \left.- \frac{1}{\tau^n_{i,j,k}}\partial_U  N_q(\Unijk)(\delta
U_{i,j,k})
-  \Delta^{n,g}_{i,j,k,q}\delta g_{i,j,k,q}
\right).
\end{split}
\end{equation*}

\clearpage

\clearpage

\clearpage

\clearpage

\clearpage

\clearpage

\clearpage

\clearpage

\clearpage


\clearpage

\clearpage

\clearpage

\clearpage

\clearpage

\clearpage

\clearpage
 
\clearpage

\clearpage

\end{document}